\documentclass[11pt]{article}
\usepackage{bbm}
\usepackage{amssymb,amsmath,amsthm}
\usepackage{graphicx}
\usepackage{flafter}
\usepackage[hidelinks]{hyperref}
\usepackage{multirow}
\usepackage{tabularx}     
\usepackage{booktabs}     

\usepackage{caption}
\usepackage{subcaption}
\usepackage{cleveref}
\allowdisplaybreaks

\usepackage{geometry}
\usepackage{color}

\usepackage[markup=underlined]{changes}
\definechangesauthor[color=blue]{ck}

\usepackage{enumitem}

\def\theequation{\thesection.\arabic{equation}}
\newtheorem{theorem}{Theorem}[section]
\newtheorem{lemma}[theorem]{Lemma}

\newtheorem{remark}[theorem]{Remark}

\theoremstyle{remark}
\theoremstyle{plain}

\def\sqr#1#2{{\vcenter{\vbox{\hrule height.#2pt
              \hbox{\vrule width.#2pt height#1pt \kern#1pt \vrule width.#2pt}
              \hrule height.#2pt}}}}

\def\Re{{\mathop{\rm Re}\,}}

\begin{document}
\title{Sharp Damping Location for Stability of a Locally Damped Wave Equation and its  Galerkin Approximation
\thanks{This work is supported by the National Natural Science Foundation of China (grants No. 12271035, 12171413, 12131008), Hebei Natural Science Foundation (grant No. A2025105012).
\medskip} }
\author{Kangsheng Liu\thanks{School of Mathematical Sciences, Zhejiang University, Hangzhou, 310058, China; Academy for Advanced Interdisciplinary Studies, Northeast Normal University, Changchun 130024, Jilin, China
(email: ksliu@zju.edu.cn).}, ~~Lihui Lu\thanks{School of Mathematical Sciences, Zhejiang University, Hangzhou, 310058, China
(email: lihui0822@zju.edu.cn).}
~~and~~Qiong Zhang\thanks{ School of Mathematics and Statistics, Beijing Institute of Technology, Beijing 100081, China; Tangshan Research Institute, Beijing Institute of Technology, Tangshan, Hebei 063000, China
(email: zhangqiong@bit.edu.cn).}
}
\date{}
\maketitle

\begin{center}
\begin{minipage}{6in}
{\bf Abstract.}
We study the stability of the variable coefficient wave equation
$u_{tt}-u_{xx}-(1+x^{2})u_{yy}+\chi_{\mathbb G}u_{t}=0$ on the rectangle
$(-1,1)\times(0,\pi)$ with Dirichlet boundary conditions, where the
viscous damping acts on a vertical strip $\mathbb G=(\eta-\delta,\eta+\delta)\times(0,\pi)$.
 We establish  a precise relationship between the location of the damping region and the rate of energy decay.
  We prove that the associated semigroup is
exponentially stable whenever the closure of the damping strip meets the line $x=0$,  i.e., $|\eta|\leq\delta$, despite the existence of rays that never enter the damped region.
When the damping
region is strictly bounded away from $x=0$, the  semigroup decays at the sharp rate $(\log(2+t))^{-1}$.
The proofs combine
separation of variables, one-dimensional resolvent estimates, a
reflection identity for the harmonic oscillator, and Gaussian
quasimodes.
  We then consider the spectral Galerkin approximation based
on the Dirichlet eigenfunctions of the constant coefficient Laplacian and show that,
although the sequence of semigroups resulting from the approximation converges to  the continuous one, it fails
to be uniformly exponentially stable.

\vskip 3mm

\vskip 4mm

 \vskip 3mm

\end{minipage}
\end{center}

\section{Introduction}

In this paper we study   a variable coefficient wave equation on the rectangular domain $\Omega = (-1,1)\times(0,\pi)$  with boundary $\partial\Omega$. Let
$$
\mathbb L u= -   u_{xx}  - (1+x^2) u_{yy}, \;\;\;\;
D(\mathbb{L}) = H^2(\Omega) \cap  H^1_0(\Omega),
$$
and consider
\begin{equation}\label{sys}
\left\{\begin{array}{ll}
u_{tt} + \mathbb{L} u +\chi_{\mathbb{G}}(x,y )u_t = 0 & \text{in } (0,\infty) \times \Omega,\\[2pt]
u = 0 & \text{on } (0,\infty) \times \partial\Omega,\\[2pt]
u(0) = u_0,\quad u_t(0) = u_1 & \text{in } \Omega,
\end{array}\right.
\end{equation}
where   the damping coefficient function $\chi_{\mathbb{G}}$ is the characteristic function of the vertical strip $
\mathbb G$  of width $2\delta$ centered at $x=\eta$:
$$
\mathbb G =  (\eta-\delta, \;\eta+\delta)\times (0,\;\pi).
$$
We assume $0 < \delta \le 1/4$ and $|\eta| \le 1-\delta$.
Since $\mathbb{L}$ is invariant under the reflection $x\to -x$,   it suffices to consider $\eta \ge 0$.

The energy associated with the system is defined by
\begin{equation}\label{energy}
E(t) = \frac12 \int_\Omega \bigl( |u_t|^2 + |u_x|^2 + (1+x^2)|u_y|^2 \bigr)\,dx\,dy.
\end{equation}
A direct computation shows that the energy satisfies the dissipation law\begin{equation}\label{derivative}
\frac{d}{dt}E(t) = - \int_\Omega \chi_{\mathbb{G}}(x,y)\,|u_t|^2 \,dx\,dy.
\end{equation}
Hence, $E(t)$ is non-increasing, and the energy is dissipated only through the viscous damping in $\mathbb{G}$.

For the classical   wave equation with localized viscous damping, the system is exponentially stable, i.e.,
$$
\Vert{}e^{tA}z\Vert{}_\mathcal H \le M e^{-\omega_0 t} \Vert{}z\Vert{}_\mathcal H,\;\;
\mbox{for some}
\; M \ge 1,\; \omega_0 > 0 \; \mbox{and all} \;  t \ge 0,
$$
 if and only if the Geometric Control Condition (GCC) holds (\cite{bardos, BG97, Ral69, RT75, L8}). This condition requires that there exists a time $T > 0$ such that every geodesic ray of length $T$ intersects the support of the damping coefficient.
 If the non-negative damping coefficient function has a positive lower bound on a nonempty subset, without satisfying the GCC, the associated semigroup is known to decay logarithmically with the   rate $(\log(2 + t))^{-1}$ (\cite{Bur98}).
 When the   domain $\Omega $ is a torus, polynomial stability can often be obtained (\cite{AL14, L10, fxy, L7, L12}).
 In particular, for a constant coefficient wave equation on a square domain with damping supported on a sub-strip, the semigroup is polynomially stable with the optimal decay rate  $(t+1)^{-2/3}$ (\cite{AL14, L6}).
See also \cite{BH07} for earlier work on   partially rectangular domains.
 General exponential and polynomial stability criteria  were developed in  \cite{L11, L3}.

In this paper, we analyze the separable, variable coefficient model \eqref{sys}.
Since   $x=0$ is the unique critical point of the coefficient function $1+x^2$,   rays near $x=0$ obey  the harmonic oscillator equation and give rise to a family of trapped rays.
We prove that the position of the damping strip determines a sharp transition in stability.
If $|\eta|\le \delta$, then the associated semigroup is exponentially stable.
In particular, even when
$|\eta|=\delta$, in which case the ray $x=0$ never enters $\mathbb{G}$,
the semigroup remains exponentially stable. On the other hand, when
$|\eta|>\delta$, we construct Gaussian quasimodes and obtain the optimal
logarithmic decay rate $(\log(2+t))^{-1}$.

The proofs combine separation of variables with one dimensional resolvent estimates. More precisely, when the damping strip reaches the line $x=0$, a rescaling  argument  reduces the high-frequency concentration problem to one for the harmonic oscillator. The reflection identity for the harmonic oscillator then forces the energy of the relevant quasimodes to be asymptotically equally distributed between the two half-lines, which is incompatible with the vanishing of the observation on the damped strip.
When the damping strip is separated from the critical line $x=0$, the localized Gaussian quasimodes have exponentially small damping, and this yields the optimal logarithmic decay rate.


We also study a   spectral Galerkin approximation  based on the Dirichlet eigenfunctions of the constant coefficient Laplacian. We refer to \cite{nicaise, ez, zua99, L14, tucsnak, L13, zua05} and references therein  for related uniform stability results.
We prove that   the corresponding finite dimensional semigroups converge strongly to the continuous semigroup generated by \eqref{sys}.   However, they fail to be uniformly exponentially stable.
The proof is based on  constructing discrete high-frequency approximate eigenvectors concentrated near $x=\pm 1$.
This loss of uniform exponential stability is due to the
fact that these eigenfunctions form only a pseudo-basis, so that the
approximating systems do not inherit the uniform stability of the
continuous one. Whether weaker types of uniform stability, such as
uniform polynomial stability, still hold is an interesting open
question.


The paper is organized as follows. Section~2 introduces the model, its well-posedness, and some preliminaries. Sections~3 and 4 prove the exponential and logarithmic stability, respectively.  Section 5 is devoted to constructing  the Galerkin approximation.

Throughout the paper, the notation $\big\langle \cdot, \cdot  \big\rangle$
denotes the inner product in $L^2(\Omega)$, and $\| \cdot \|$ denotes
the associated norm. The letter $C$ stands for a generic positive constant
that may vary from line to line.

\section{Preliminaries}
\setcounter{equation}{0} \setcounter{theorem}{0}

In this section, we reformulate system \eqref{sys} as an abstract Cauchy
problem and present some preliminaries.
We define the bounded linear operator $\mathbb B$ on $ L^2(\Omega)$ by $\mathbb B u = \chi_{\mathbb{G}}u$. System \eqref{sys} can then be written in the following   form
\begin{equation} \label{21}
      u_{tt}+ \mathbb L u  + \mathbb B u_t = 0,\;\; u(0)=u_0,\;\; u_t(0)=u_1.
\end{equation}
We introduce the energy   space
$$ \mathcal{H} = H_0^1(\Omega) \times L^2(\Omega), $$
equipped with the   inner product
\begin{equation}
    \big\langle U_1, U_2 \big\rangle_{\mathcal{H}} = \big\langle \mathbb L^{1\over2} u_1,\mathbb L^{1\over2} u_2 \big\rangle + \big\langle v_1, v_2 \big\rangle,
\end{equation}
for all $U_1 = (u_1, v_1)^T$ and $U_2 = (u_2, v_2)^T$ in $\mathcal{H}$.

With the state vector $U(t) = (u(t), u_t(t))^T$, system \eqref{sys} can be
rewritten as the abstract first-order evolution equation on $\mathcal{H}$:
\begin{equation} \label{abs}
    \frac{\mathrm{d}}{\mathrm{d}t} U(t) = A U(t), \quad U(0) = (u_0, u_1)^T \in \mathcal{H}.
\end{equation}
Here, the unbounded linear operator $A: \mathcal{D}(A) \subset \mathcal{H} \to \mathcal{H}$ is defined by
\begin{equation}
    A \begin{pmatrix} u \\ v \end{pmatrix} = \begin{pmatrix} v \\ -\mathbb{L}u - \mathbb B v \end{pmatrix},
\end{equation}
with domain $D(A) = D(\mathbb{L}) \times H_0^1(\Omega)$.

For any $U = (u,v)^T \in D( A )$, a direct computation yields
\begin{equation}
    \Re \big\langle AU, U \big\rangle_{\mathcal{H}} =   - \int_{\mathbb{G}} |v|^2 \, \mathrm{d}x\mathrm{d}y \le 0.
\end{equation}
This shows that $A$ is dissipative. Since $A$ is closed and densely
defined, and $0 \in \rho(A)$,  the Lumer--Phillips theorem
 implies that $A$ generates a $C_0$-semigroup of contractions
$\{e^{tA}\}_{t \ge 0}$ on the energy space $\mathcal{H}$ (\cite{engel,L3}).

The strong stability of system \eqref{sys} was established in \cite{chenliuliu, L3}.
\begin{lemma}\label{lemma-stro}
  $i\mathbb{R} \subset \rho(A)$. Consequently, the semigroup $\{e^{tA}\}_{t\ge0}$ is strongly stable on $\mathcal{H}$.
\end{lemma}

In this paper, we investigate how the stability of system  \eqref{sys}  depends on the location of the damping strip $\mathbb{G}$  and, in particular,  on whether it intersects the $y$-axis.
We show that even  if  there are rays that never enter the damping region, the system  remains  exponentially stable  provided that  $\mathbb{G}$ reaches the $y$-axis.
When $\mathbb{G}$  is bounded away from  the $y$-axis, however, the system  is merely  logarithmically stable, and  this decay rate is  optimal.

The proofs of exponential and logarithmic stability rely  on frequency-domain techniques, which reduce the stability of the semigroup to resolvent estimates on the imaginary axis. The following lemma  summarizes these classical criteria.

\begin{lemma}\label{lemma-bt}
 Let $A$ be the generator of a uniformly bounded $C_0$-semigroup $\{e^{tA}\}_{t\ge0}$ on a Hilbert space $H$.
\begin{enumerate}
\item[\em (i)] (\cite{lyubich_vu}) If $i\mathbb{R} \subset \rho(A)$, then $\{e^{tA}\}_{t\ge0}$ is strongly stable.

\item[\em (ii)]   (\cite{gearhart, huang, pruss})
The semigroup is exponentially stable if and only if $i\mathbb{R} \subset \rho(A)$ and$$\limsup_{\vert{}\omega\vert{}\to\infty} \Vert{}(i\omega I - A)^{-1}\Vert{}_{\mathcal{L}(H)} < \infty.$$

\item[\em (iii)] (\cite{bt})
 For a given $\alpha > 0$, the semigroup is polynomially stable of order $\alpha$ if and only if $i\mathbb{R} \subset \rho(A)$ and$$\limsup_{\vert{}\omega\vert{}\to\infty} \vert{}\omega\vert{}^{-1/\alpha} \Vert{}(i\omega I - A)^{-1}\Vert{}_{\mathcal{L}(H)} < \infty.$$

\item[\em (iv)] (\cite{batty_duyckaerts})
For a given $\alpha > 0$, the semigroup is logarithmically stable of order $\alpha$ if and only if $i\mathbb{R} \subset \rho(A)$ and there exists a constant $c > 0$ such that$$\limsup_{\vert{}\omega\vert{}\to\infty} e^{-c\vert{}\omega\vert{}^{1/\alpha}} \Vert{}(i\omega I - A)^{-1}\Vert{}_{\mathcal{L}(H)} < \infty.$$

\end{enumerate}

\end{lemma}

An equivalent condition for the exponential stability of second-order systems was established in \cite{L3} based on Lemma \ref{lemma-bt} (ii). Applying this result to system \eqref{21} yields the following lemma.

\begin{lemma} \label{lemma-lk}
 System \eqref{21} is exponentially stable if and only if  there exists a constant  $C_\delta > 0$ such that for all $\omega \in \mathbb{R}$ and all $u \in H^2(\Omega)\cap H_0^1(\Omega)$,
\begin{equation} \label{freq}
    \|(\omega^2 - \mathbb{L})u\|  + |\omega|\|u\|_{L^2(\mathbb{G})} \ge C_\delta \|\mathbb{L}^{{1\over2}}u\| .
\end{equation}
 \end{lemma}

  A natural question is whether finite-dimensional approximations of \eqref{sys} converge to the continuous system and inherit its stability properties.
 We  introduce the following orthonormal basis of $L^2(\Omega)$.
  \begin{align} \label{basis}
 \widetilde{e}_m(x) := \sin\frac{m\pi(x+1)}{2}, \quad e_n(y) :=   \sqrt{\frac{2}{\pi}}  \sin(ny), \quad \phi_{m,n}(x,y) := \widetilde{e}_m(x)e_n(y), \;\; (x,y)\in \Omega.
 \end{align}
It is known that
 $\{\phi_{m,n}\}_{m,n\ge1}$
forms an orthonormal basis of $L^2(\Omega)$ consisting of eigenfunctions of
  $$\mathbb L_0  = - \partial_{xx}    - \partial_{yy}, \;\;\;\;
D(\mathbb{L}_0) = H^2(\Omega) \cap  H^1_0(\Omega),  $$ corresponding to  the  eigenvalues $ \{\lambda_{m,n}\}_{m,n\ge1} = \big\{\big( {m\pi / 2}\big)^2   + n^2\big\}_{m,n\ge1}$.
The approximation  spaces are constructed from   eigenfunctions  $\{\phi_{m,n}\}_{m,n\ge1}$,
rather than those of the variable coefficient operator $\mathbb L$, which does not have explicit eigenfunctions. For this reason, we refer to this basis as a pseudo-eigenbasis of \eqref{sys}.

  \begin{remark}\label{remark-ll}
For any \(z \in H^2(\Omega) \cap H_0^1(\Omega)\), one has
\begin{align*}
\|\mathbb L z\|^2  = \sum\limits_{n=1}^\infty\big(\|z_n''\|^2_{L^2(-1,1)} + n^4\|(1+x^2)z_n\|^2_{L^2(-1,1)}
  + 2n^2\|\sqrt{1+x^2}z_n'\|^2_{L^2(-1,1)}   - 2n^2\|z_n\|^2_{L^2(-1,1)}\big),
\end{align*}
where $z_n \in H^2(-1,1) \cap H_0^1(-1,1)$ are the Fourier coefficients of $z$  with respect to  the orthonormal basis $e_n(y)$ of $L^2(0,\pi)$.
Since  \(\|z_n\|_{L^2(-1,1)} \leq (2/\pi)\|z_n'\|_{L^2(-1,1)}\) for each $n$,
\[
\|\mathbb L z\|^2 \geq \sum\limits_{n=1}^\infty\big(\|z_n''\|^2_{L^2(-1,1)} + n^4\|(1+x^2)z_n\|^2_{L^2(-1,1)}
  + 2(1- 4/\pi^2)n^2    \|z_n'\|^2_{L^2(-1,1)} \big).
\]
Consequently, there exists a constant $C > 0$ such that
\[
\|\mathbb L_0 u\| \leq C\|\mathbb Lu\|, \qquad \|\mathbb Lu\| \leq 2\|\mathbb L_0 u\|.
\]
 \end{remark}

 By constructing appropriate approximation spaces, we show that the corresponding finite dimensional semigroups converge to the continuous semigroup $\{e^{tA}\}_{t\ge0}$. However, the family of Galerkin approximations  fails to  inherit the  exponential stability of the continuous system (Theorem~\ref{th-ne}).  The analysis relies on  the following lemma  (\cite{lks-l, lzy-zhsm}).

\begin{lemma}  \label{lemma-g2} 
Let \( \{e^{tA_n}\}_{t\ge 0 } \)  be a sequence of contraction \( C_0 \)-semigroups on Hilbert spaces \( H_n, \; n=1,2,\cdots, \) and let  \( A_n \) be the corresponding infinitesimal generators. Then the following  are equivalent:
\begin{enumerate}
 \item[{\rm (i)}]  The family of semigroups  \(\{e^{tA_n}\}_{t\ge 0 }  \)   is uniformly exponentially stable, i.e.,
     there exist  \( M \ge1 \;\; \omega_0>0 \) (independent of \( n \)) such that
\[
\sup_{n\in {\mathbb N}}\|e^{tA_n}\|_{\mathcal{L}(H_n)} \le Me^{-\omega_0 t}, \quad t \ge 0.
\]

\item[{\rm (ii)}] $i\mathbb{R} \subset \rho(A_n)$ for all $n \in \mathbb{N}$, and the family of resolvents   is uniformly  bounded on the imaginary axis, i.e.,
\[
\sup\limits_{\omega  \in  \mathbb{R}, \, n \in \mathbb{N} }   \| (i \omega  - A_n)^{-1} \|_{\mathcal{L}(H_n)}  < \infty.
\]
\end{enumerate}
\end{lemma}


\section{Exponential stability when $0\le \eta \le \delta$}
\setcounter{equation}{0} \setcounter{theorem}{0}

In this section, we prove exponential stability of system \eqref{sys} when  $0\le \eta \le \delta$.
In this case, the closure of the damping strip always intersects the $y$-axis. In particular, the rightmost admissible damping strip is
$
(0,2\delta)\times(0,\pi)
$.
The main result in this section reads as follows.

\begin{theorem} \label{th2.1} Assume that $0 < \delta \le 1/4$ and   $0 \le \eta \le \delta$. Then the semigroup $\{e^{tA}\}_{t\ge0}$ is exponentially stable.
\end{theorem}

By Lemma \ref{lemma-lk},  it remains to  prove the   estimate \eqref{freq}. Since $G_\delta := (0,\delta)\times(0,\pi) \subset \mathbb{G}$, it is sufficient to establish \eqref{freq} with the damping region $\mathbb{G}$ replaced by $G_\delta$.
Suppose that  \eqref{freq} fails. Then there exist a sequence of real numbers $\omega_k \in \mathbb{R}$ and a sequence of functions $u_k \in H^2(\Omega)\cap H_0^1(\Omega)$ such that
as $k \to \infty$,
\begin{equation}\label{e1}
    \|\mathbb{L}^{{1\over2}}u_k\|  = 1, \qquad \|F_k \|  = \|(\omega_k^2 - \mathbb{L})u_k\|   = o(1), \qquad \|\omega_k u_k\|_{L^2(G_\delta)} = o(1).
\end{equation}

We first show that $|\omega_k| \to \infty$. Indeed, if the sequence $\{\omega_k\}$ were bounded, we could extract a subsequence (still denoted by $\omega_k$) converging to some $\omega_0 \in \mathbb{R}$. Since $i\omega_0 \in \rho(A)$ by Lemma \ref{lemma-stro}, the resolvent $(i\omega - A)^{-1}$ is uniformly bounded in a neighborhood of $i\omega_0$.
 Since \(\mathbb Lu_k = \omega_k^2 u_k - F_k\) is bounded in $L^2(\Omega)$, taking a subsequence \(u_k \to u\) in \(H_0^1\), we have \(\mathbb Lu = \omega_0^2 u\), \(\|\mathbb L^{1/2}u\| = 1\), and \(\omega_0 u = 0\) in \(G_\delta\). If \(\omega_0 = 0\), positivity of $\mathbb L$ yields a contradiction. Otherwise, the uniqueness for the Cauchy problem of the ODE satisfied by each Fourier mode yields $u \equiv 0$, a contradiction.

Since the coefficient $1+x^2$ depends only on $x$ and the domain is a rectangle, we can expand $u_k$ and $F_k$ in the orthonormal basis $\{e_n(y)\}_{n=1}^{\infty}$ of $L^2(0,\pi)$ as follows:
\begin{equation}\label{e2}
    u_k(x,y) = \sum_{n \ge 1} u_{k,n}(x)e_n(y), \quad F_k(x,y) = \sum_{n \ge 1} f_{k,n}(x)e_n(y),
\end{equation}
where $u_{k,n}(x)$ and $f_{k,n}(x)$ are the Fourier coefficients and $e_n(y)$ is defined in \eqref{basis}.
 Define the energy of the $n$-th mode by
\begin{equation}
    E_{k,n} := \int_{-1}^{1} \left( |u_{k,n}'(x)|^2 + n^2(1+x^2)|u_{k,n}(x)|^2 \right) \mathrm{d}x,
\end{equation}
and the combined residual and damping term by
\begin{equation}
    R_{k,n} := \|f_{k,n}\|_{L^2(-1,1)}^2 + |\omega_k|^2 \|u_{k,n}\|_{L^2(0,\delta)}^2.
\end{equation}
It then follows from \eqref{e1} that
\begin{equation}\label{sum-eq}
    \sum_{n \ge 1} E_{k,n} = 1, \qquad \sum_{n \ge 1} R_{k,n} = \|F_k\| ^2 + \|\omega_k u_k\|_{L^2(G_\delta)}^2 =: \epsilon_k^2 = o(1).
\end{equation}
Since   $\sum_{n \ge 1} E_{k,n} = 1$, the sum $\sum_{n \ge 1} R_{k,n}$ is a convex combination of $R_{k,n} / E_{k,n}$ with weights $E_{k,n}>0$.
Therefore, there exists an index \(n_k\) with $E_{k,n_k} > 0$ such that
\begin{equation} \label{e2.10}
\frac{R_{k,n_k}}{E_{k,n_k}}
\le
\epsilon_k^2 =o(1).
\end{equation}
We now define the energy normalized one dimensional sequences associated with this mode:
\begin{equation}
    \tilde{u}_k(x) := \frac{u_{k,n_k}(x)}{\sqrt{E_{k,n_k}}}, \qquad \tilde{f}_k(x) := \frac{f_{k,n_k}(x)}{\sqrt{E_{k,n_k}}}.
\end{equation}
By \eqref{e2.10},   $\|\tilde{f}_k\|_{L^2(-1,1)}^2 + |\omega_k|^2 \|\tilde{u}_k\|_{L^2(0,\delta)}^2 \le \epsilon_k^2 = o(1)$. In what follows, we drop the tildes and write $u_k$, $f_k$, and $n_k$ for notational simplicity.


The following proof is divided into several lemmas.

\begin{lemma}\label{lemma-e1}
The sequences $u_k,\; f_k,\; \omega_k,\; n_k$  satisfy
\begin{align} \label{norm}
&u_{k}(\pm1)=0,\qquad
\int_{-1}^1\big(|u_{k}'|^2+n_k^2(1+x^2)|u_{k}|^2\big)\,dx=1,
\end{align}
and as $k\to \infty,$
\begin{align} \label{resi}
&\|f_{k}\|_{L^2(-1,1)}=o(1),\qquad f_{k}= u_{k}''+q_ku_{k},\qquad q_k(x) := \omega_k^2 - n_k^2(1+x^2), \\
\label{damp}
&|\omega_k|\,\|u_{k}\|_{L^2(0,\delta)}=o(1),
\\ \label{e3}&
\omega_k^2 \|u_{k}\|_{L^2(-1,1)}^2 = 1 + o(1),
\\ \label{ee} &
{n_k  \over |\omega_k|} \le 1 + o(1).
 \end{align}
 \end{lemma}
 \begin{proof}
 Properties \eqref{norm}-\eqref{damp} follow directly from \eqref{e1}-\eqref{e2}.
Taking the inner product of \eqref{resi} with $u_{k}$ and integrating by parts, we obtain
  \begin{align}
  \int_{-1}^1\big(\omega_k^2 |u_{k}|^2 -|u_{k}'|^2 - n_k^2(1+x^2)|u_{k}|^2\big)\,dx= o(1).
  \end{align}
  Combining this identity with \eqref{norm} yields \eqref{e3}.
Moreover, by \eqref{norm},
$
n_k^2\|u_k\|_{L^2(-1,1)}^2 \le 1.
$
Together with \eqref{e3}, this yields \eqref{ee}.
\end{proof}

By the uniform boundedness in \eqref{ee}, we can extract a subsequence (still indexed by $k$) such that
\begin{equation}\label{e27_theta1}
    \lim_{k\rightarrow\infty} \frac{n_k}{|\omega_k|} = \theta \in [0,1].
\end{equation}


\begin{lemma}\label{lemma-e2}
Assume $\theta < 1$. Then as $ k \to \infty$,
\begin{equation}
\label{thetasmall}
\int^{1 }_{-1} \left(|u_{k}'|^2 +  \omega_k^2 |u_{k}|^2 \right) \,dx = o(1).
\end{equation}
\end{lemma}
\begin{proof}

 Let $0 < \rho < \delta$ and let $\zeta \in C_0^\infty(0,\rho)$ be a smooth cut-off function satisfying $0 \le \zeta \le 1$ and $\zeta \equiv 1$ on $[\varepsilon, \rho - \varepsilon]$ for some $0 < \varepsilon < \rho/2$. Multiplying \eqref{resi} by $\zeta^2 \bar{u}_k$ and integrating over $(0,\rho)$,   we obtain
\begin{equation*}
    \int_0^\rho \zeta^2 |u_k'|^2 \,dx = \int_0^\rho q_k \zeta^2 |u_k|^2 \,dx - 2   \int_0^\rho \zeta \zeta' u_k' \bar{u}_k \,dx -  \int_0^\rho f_k \zeta^2 \bar{u}_k \,dx.
\end{equation*}
Applying the Cauchy-Schwarz and Young  inequalities to the middle term, it follows that
\begin{equation} \label{e51}
    \frac{1}{2} \int_0^\rho \zeta^2 |u_k'|^2 \,dx \le \int_0^\rho \left( q_k \zeta^2 + 2 |\zeta'|^2 \right) |u_k|^2 \,dx + \int_0^\rho \zeta^2 |f_k| |u_k| \,dx.
\end{equation}
Recall   that $q_k(x) \le \omega_k^2$. Since $\operatorname{supp}(\zeta) \subset (0,\rho) \subset (0,\delta)$, the damping condition \eqref{damp} implies that  every term on the right-hand side of \eqref{e51} converges to zero. Thus, we obtain
\begin{equation} \label{e5}
    \int_\varepsilon^{\rho-\varepsilon} |u_k'|^2 \,dx \le \int_0^\rho \zeta^2 |u_k'|^2 \,dx \to 0 \quad \text{as } k \to \infty.
\end{equation}

  For any $C^1$ function $h$ on an interval $I \subset [-1,1]$, multiplying \eqref{resi} by $h \bar{u}_k'$, taking the real part, and integrating by parts over $I$, we obtain
\begin{equation} \label{e_rellich}
\begin{aligned}
2\operatorname{Re} \int_{I} f_k h \bar{u}_k' \,dx &= \left[ h |u_k'|^2 + h q_k |u_k|^2 \right]_{\partial I} - \int_{I} \left( h'|u_k'|^2 + (h q_k)'|u_k|^2 \right) \,dx.
\end{aligned}
\end{equation}

Let $\xi \in C^1([0,\rho])$ be a cut-off function satisfying $0 \le \xi \le 1$ and
\begin{equation*}
    \xi \equiv 1 \;\; \mbox{on} \;\; [0, \rho/4], \qquad \xi \equiv 0 \;\; \mbox{on} \;\; [3\rho/4, \rho].
\end{equation*}
By construction, $\operatorname{supp}(\xi') \subset [\rho/4, 3\rho/4]$. Choosing  $\varepsilon$   small enough such that $\varepsilon \le \rho/4$, we can apply \eqref{e_rellich} on   $I = (0,\rho)$ with $h = \xi$. The boundary term at $x=\rho$ vanishes since $\xi(\rho)=0$. Using the uniform boundedness of $n_k/|\omega_k|$, along with \eqref{damp} and \eqref{e5}, we deduce
\begin{equation} \label{e52}
    |u_k'(0)|^2 + q_k(0) |u_k(0)|^2 = o(1) \quad \text{as } k \to \infty.
\end{equation}
 Define
$$\mathcal V_{k}(x) := |u_{k}'(x)|^2 + q_k(x)|u_{k}(x)|^2.$$
  Taking $I = (0,1)$ and the multiplier $h(x) = x-2$ in \eqref{e_rellich}, we have
\begin{equation}
    -|u_k'(1)|^2 + 2\mathcal{V}_k(0) - \int_{0}^1 \left( |u_k'|^2 + ((x-2)q_k)'|u_k|^2 \right) \,dx = o(1).
\end{equation}
Similarly, taking $I = (-1,0)$ and $h(x) = x+2$ in \eqref{e_rellich} yields
\begin{equation}
    -|u_k'(-1)|^2 + 2\mathcal{V}_k(0) - \int_{-1}^0 \left( |u_k'|^2 + ((x+2)q_k)'|u_k|^2 \right) \,dx = o(1).
\end{equation}
 Summing these two identities,   using $ ((x\mp 2)q_k)' = \omega_k^2 - n_k^2 - 3n_k^2 x^2 \pm 4n_k^2 x$  and \eqref{e52}, we obtain
\begin{equation}
\label{e260}
  |u_{k}'(-1)|^2+  |u_{k}'(1)|^2    +  \int^{1 }_{-1} \left[ |u_{k}'|^2 +  (\omega_k^2   + n_k^2 (  - 3  x^2 + 4|x|-1))|u_{k}|^2 \right] \,dx = o(1).
\end{equation}
  Note that $ - 3  x^2 + 4|x|-1 \in [-1, \; 1/3] $ for $x\in [-1,1]$. Furthfermore,
   it follows from $\theta <1 $ that $\theta^2 < {{1+\theta^2}\over{2}} <1,$  and   consequently,
   $$
   n_k^2 \le  {{1+\theta^2}\over{2}} \omega_k^2, \;\; \;  \mbox{ for   }  \; k \;  \mbox{  large enough. }  \;
   $$
   Combining this with \eqref{e260} yields
   \begin{equation}
\label{e261}
  0\le  |u_{k}'(-1)|^2+  |u_{k}'(1)|^2    +  \int^{1 }_{-1} \left( |u_{k}'|^2 +  {{1-\theta^2}\over{2}}|u_{k}|^2 \right) \,dx = o(1).
\end{equation}
  Consequently, we obtain \eqref{thetasmall}.
\end{proof}

From \eqref{thetasmall}, we arrive at a contradiction to \eqref{e3} in the case $\theta<1$.
\medskip

We now consider the remaining case where $\theta = 1$.
In the following two lemmas, we show that the damping condition forces the right half-energy to vanish,  while the harmonic oscillator limit yields an asymptotic equality of the
energy on the two half-lines. These two conclusions contradict each other, and hence the case $\theta = 1$  cannot occur.

 \begin{lemma}
 \label{lemma-ykn} Assume $\theta = 1$.
For the sequence $n_k$ introduced in \eqref{e2.10}, set
\begin{align}
y_{k}(X):=
n_k^{-{1\over4}} |\omega_k| u_{k}  (n_k^{-{1\over2}}  X ),
\qquad \;  \mbox{  for }  \;
|X|<n_k^{1\over2},\;\; k\in {\mathbb N}.
 \end{align}
Then,  as  $k \to \infty$,  the following two  conclusions hold simultaneously, which is impossible.
\begin{align} \label{e291}
& \int_{0}^{n_k^{1\over2}}
|y_{k}|^2 dX
=o(1),
 \\ \label{e2911}
& \int_{0}^{n_k^{1\over2}}
|y_{k}|^2 dX
=
{1\over2}+o(1), \;\; \;\;
\int_{-n_k^{1\over2}}^{0}
|y_{k}|^2 dX
=
{1\over2}+o(1).
\end{align}
 \end{lemma}

 \begin{proof} {\bf (i)}
By \eqref{norm} and \eqref{resi}, $y_{k}$ satisfies
the boundary conditions
$
y_{k}(\pm n_k^{1\over2})=0,
$
and
\begin{align} \label{e20}
\left(
H  -
\Lambda_k
\right)
y_{k}(X)=
- |\omega_k| n_k^{-{5\over4}} f_{ k}(n_k^{-{1\over2}} X),\;\;\;\;  \forall \; X\in (-n_k^{1\over2}, n_k^{1\over2}),
\end{align}
where
\begin{align}\label{hh}
 H := -\partial_X^2 +
X^2, \;\;\;\; \Lambda_k :=
\frac{\omega_k^2-n_k^2}{n_k},
 \end{align}
 and
 $$
\big\||\omega_k| n_k^{-{5\over4}} f_{ k}(n_k^{-{1\over2}} X)\big\|_{L^2(-n_k^{1\over2}, n_k^{1\over2})} = |\omega_k| n_k^{-1} \|f_k\|_{L^2(-1,1)} = o(1).
 $$

It follows from \eqref{damp} that
 \begin{align} \label{edd}
 \int_{0}^{\delta n_k^{1\over2}}
|y_{k}(X)|^2 dX =o(1).
 \end{align}
Moreover, \eqref{e3} gives
\begin{align}\label{ynorm}
 \| y_k \|_{L^2(-n_k^{1\over2}, n_k^{1\over2})} = 1+o(1).
\end{align}
Taking the $L^2(-n_k^{1\over2}, n_k^{1\over2})$ inner product of \eqref{e20} with $y_k$ and using \eqref{e27_theta1}, \eqref{ynorm}, we obtain
 \begin{align} \label{ed}\begin{split}
& \int_{-n_k^{1\over2}}^{n_k^{1\over2}}
 (
|y_{k}' (X)|^2
+
X^2|y_{k} (X)|^2
 ) dX \\ =\; &
\Lambda_k\int_{-n_k^{1\over2}}^{n_k^{1\over2}}
   |y_{k} (X)|^2 dX   - |\omega_k| n_k^{-{5\over4}} \int_{-n_k^{1\over2}}^{n_k^{1\over2}}
  f_{k}( n_k^{-{1\over2}} X ) y_{k} (X) dX
\\ =\; &
\Lambda_k\int_{-n_k^{1\over2}}^{n_k^{1\over2}}
   |y_{k} (X) |^2 dX
+
o(1).
\end{split}  \end{align}
By \eqref{e27_theta1},
 \begin{align}\label{ed1} \omega_k^2-n_k^2 = o(n_k^2), \;  \mbox{ and   }  \; \Lambda_k =    o(n_k).
 \end{align}
 Therefore, from \eqref{ynorm}, \eqref{ed} and \eqref{ed1},
  \begin{align}\label{ed0}\begin{split}
  \int_{\delta n_k^{1\over2}}^{n_k^{1\over2}}
|y_{k}|^2dX  & \le {1\over \delta^2 n_k}\int_{\delta n_k^{1\over2}}^{n_k^{1\over2}} X^2
|y_{k}|^2dX  \\ \noalign{\medskip}
  & \le
\frac{\Lambda_k \| y_k \|^2_{L^2(-n_k^{1\over2}, n_k^{1\over2})} +o(1)}
{\delta^2 n_k} \to 0.
\end{split}\end{align}
Combining \eqref{edd} and \eqref{ed0}, we arrive at \eqref{e291}.

 {\bf (ii)}  Choose a cutoff function $\chi \in C_c^\infty(-1,1)$ such that $\chi \equiv 1$ on $[-1/2, 1/2]$, and define
\begin{align} \label{y37}
\chi_k(X) = \chi\big(n_k^{-{1\over2}}X\big) \qquad \mbox{for } X \in \big[-n_k^{1\over2},\; n_k^{1\over2}\big].
\end{align}
Then, we have
\begin{align} \label{ch39}
\chi_k(X) \equiv 1 \qquad \mbox{for } X \in \big[-n_k^{1\over2} / 2,\; n_k^{1\over2} / 2 \big].
\end{align}
As in the proof of  \eqref{ed0}, since $X^2 \ge n_k/4$ on the support of $1-\chi_k$, we   deduce from \eqref{ynorm}-\eqref{ed1} that
\begin{align}\label{y39}
\int_{ \{n_k^{1\over2}/2 \le |X| \le n_k^{1\over2} \} } |y_k|^2 dX
\le \frac{4}{n_k} \int_{ \{n_k^{1\over2}/2 \le |X| \le n_k^{1\over2}\} } X^2 |y_k|^2 dX \to 0.
\end{align}
Furthermore, by \eqref{ynorm}, \eqref{ed}, \eqref{ed1}, and \eqref{y37}, we estimate the derivative term:
\begin{align} \label{y40}
\begin{split}
\|\chi_k' y_{k}'\|_{L^2(-n_k^{1\over2}, n_k^{1\over2})}
& \le \|\chi_k'\|_{L^\infty(-n_k^{1\over2}, n_k^{1\over2})} \|y_{k}'\|_{L^2(-n_k^{1\over2}, n_k^{1\over2})} \\
& \le C n_k^{-{1\over2}} \big( \Lambda_k \|y_k\|_{L^2(-n_k^{1\over2}, n_k^{1\over2})}^2 + o(1) \big)^{1\over2} \\
& = O\Big( \big(\Lambda_k / n_k\big)^{1\over2} \Big) = o(1).
\end{split}
\end{align}
Extend $\widetilde y_{k} := \chi_k y_{k}$ by zero to all of $\mathbb{R}$, i.e.,
$$
\widetilde y_{k}(X) =
\begin{cases}
\chi_k(X) y_{k}(X), & X \in [-n_k^{1\over2},\; n_k^{1\over2}], \\
0, & X \in \mathbb{R} \setminus [-n_k^{1\over2},\; n_k^{1\over2}].
\end{cases}
$$
Then, from \eqref{ynorm} and \eqref{y39}, it follows that
\begin{align} \label{y41}
\|\widetilde y_{k}\|_{L^2(\mathbb{R})} = 1+o(1).
\end{align}
Moreover,   note that
\begin{align}\label{y42}
\big(H - \Lambda_k\big)\widetilde y_{k} = \chi_k \big(H - \Lambda_k\big)y_{k} - 2\chi_k' y_{k}' - \chi_k'' y_{k}, \qquad X \in [-n_k^{1\over2},\; n_k^{1\over2}].
\end{align}
By \eqref{e20}, \eqref{ch39}--\eqref{y40}, and the fact that $\|\chi_k''\|_{L^\infty({\mathbb R})} = O(n_k^{-1})$, we obtain
\begin{align} \label{d45}
\big\|\big(H - \Lambda_k\big)\widetilde y_{k}\big\|_{L^2(\mathbb{R})} =o(1).
\end{align}
A direct computation shows that
\begin{align}
(i \partial_\tau -H) (e^{-i\tau H}\widetilde y_{k} - e^{-i\tau\Lambda_k}\widetilde y_{k}) = (H-\Lambda_k ) e^{-i\tau\Lambda_k}\widetilde y_{k}.
 \end{align}
Hence, by Duhamel's formula, the unitarity of  the group $e^{-i\tau H}$, and \eqref{d45}, we have
\begin{align} \label{d47}
 \|
e^{-i\tau H}\widetilde y_{k}
-
e^{-i\tau\Lambda_k}\widetilde y_{k}
 \|_{L^2({\mathbb R})} & = \left\| - i \int_0^\tau e^{-i(\tau-s)H} (H-\Lambda_k)e^{-is\Lambda_k}\widetilde y_{k}  ds \right\|_{L^2({\mathbb R})} \\ &
\le
|\tau| \|
(H-\Lambda_k)\widetilde y_{k}
\|_{L^2({\mathbb R})}
= o(1).\nonumber
\end{align}

The harmonic oscillator $H = -\partial_X^2 + X^2$ on $L^2(\mathbb{R})$ has purely discrete spectrum $\{ 2j+1 \,:\, j \in \mathbb{N} \cup \{0\} \}$, with corresponding eigenfunctions $\xi_j(X) = h_j(X) e^{-X^2/2}$, where $h_j$ is the Hermite polynomial satisfying $h_j(-X) = (-1)^j h_j(X)$; see, e.g., \cite{thang}.

 Define the reflection operator $R f(X) = f(-X)$ for $X \in \mathbb{R}$. Then the Hermite eigenfunctions satisfy
  \begin{align}
  R\xi_j = (-1)^j \xi_j.
   \end{align}
   Moreover, it is clear that 
\begin{align} \label{hre}
e^{-i\pi H/2} \xi_j=e^{-i\pi (2j+1)/2} \xi_j=
e^{-i\pi/2}R\xi_j,\;\;\;\; \forall j \in \mathbb{N} \cup \{0\}.
\end{align}
Since the Hermite eigenfunctions form an orthogonal basis for $L^2(\mathbb{R})$,
 setting $\tau = \pi/2$ in \eqref{d47} and applying \eqref{hre}, we obtain
\begin{align}
\big\| e^{-i\pi/2} R \widetilde y_{k} - e^{-i\pi\Lambda_k/2}\widetilde y_{k} \big\|_{L^2(\mathbb{R})} =o(1).
\end{align}
Consequently,   with  $c_k := e^{-i\pi(\Lambda_k - 1)/2}$ so that $|c_k| = 1$,
\begin{align}
\|
R\widetilde y_{k}-
c_k\widetilde y_{k}\|_{L^2({\mathbb R})}
\to 0.
\end{align}
  Combining this  with \eqref{y41} yields
\begin{align} \label{e511}
\int_{0}^{\infty}
|\widetilde  y_{k}|^2 dX
=
{1\over2} +o(1), \;\; \;\;
\int_{-\infty}^{0}
|\widetilde  y_{k}|^2 dX
=
{1\over2} +o(1).
\end{align}
Then,
 by the definition of $\widetilde  y_{k}$ and  \eqref{y39}, \eqref{e511}, we obtain   \eqref{e2911}. The proof is complete.
 \end{proof}

 From Lemma \ref{lemma-ykn}, we conclude that when  $\theta = 1$,  the following two incompatible conclusions hold   as $ k \to \infty$,
\begin{align} \label{e29111}
  &  \int_{0}^{1} \omega_k^2 |u_k|^2 \,dx = \frac{1}{2}+o(1), \qquad \int_{-1}^{0} \omega_k^2 |u_k|^2 \,dx = \frac{1}{2}+o(1),
\\ & \label{e29121}
    \int_{0}^{1} \omega_k^2 |u_k|^2 \,dx = o(1) .
\end{align}
In summary, the desired contradiction follows from \eqref{e29111} and \eqref{e29121}, and the proof of Theorem \ref{th2.1}  is complete.

\begin{remark}
In the above proof, the assumption \(\eta \le \delta\) is used  in the inclusion \(G_\delta = (0,\delta)\times(0,\pi) \subset \mathbb G\).
 \end{remark}

\section{Optimal logarithmic rate when $\delta < \eta \le 1-\delta$}
\setcounter{equation}{0} \setcounter{theorem}{0}

In this section we consider the case $\delta<\eta\le1-\delta$   so that the damping strip $\mathbb{G}$ lies strictly to the right of the $y$-axis. We  show that the system \eqref{sys} is logarithmically stable with the optimal decay rate.

\begin{theorem} \label{th2.2} Assume that $0 < \delta \le 1/4$ and $\delta < \eta \le 1-\delta$. Then the semigroup $\{e^{tA}\}_{t\ge0}$  is logarithmically stable with the decay estimate
$$\|e^{tA}U\|_{\mathcal H} \le \frac{C}{\log(2+t)}\|AU\|_{\mathcal H},
\qquad t\ge0,\; U \in D(A). $$
 \end{theorem}

First, we show  that the semigroup $e^{tA}$ is logarithmically stable. From  Lemma \ref{lemma-bt}, it suffices to establish the following resolvent estimate.

\begin{lemma}\label{lemma-log}
For every fixed $\eta$ with $\delta<\eta\le1-\delta$, there exist constants $C_\eta, \; \tilde{C}_\eta>0$ such that
\begin{align}\label{lr}
\|(i\omega I-A)^{-1}\|_{\mathcal L (\mathcal H)}
\le
C_\eta e^{\tilde{C}_\eta|\omega| },
\qquad
\forall \; \omega \in\mathbb{R}.
 \end{align}
\end{lemma}

\begin{proof}For $\omega \in {\mathbb R}$ and $(f,g)^\top\in \mathcal H$, consider the equation
\begin{equation}
(i\omega I - A)U =F = (f,g)^\top,\;\;\;\;  U =  (u,v)^\top \in D(A) .
 \end{equation}
It follows that
\begin{equation}
\label{lyu}
 v= i\omega u -  f, \;\;\;\;
\omega^2 u - \mathbb{L}u   - i\omega \chi_{\mathbb{G}}(x,y)u = F :=-   g - (i\omega + \chi_{\mathbb{G}}(x,y))f.
\end{equation}
Expanding $u, \;  v, \;  F$ as in  \eqref{e2}, we find that the Fourier coefficients $u_n,\;  F_n$ satisfy
\begin{equation}
\label{le}
 u_{n}'' + q_n(x) u_{n} - i\omega \chi_{I_{\eta, \delta}}(x)u_n = F_{n},
\end{equation}
 where $I_{\eta, \delta} :=  (\eta- \delta, \; \eta+\delta),\; q_n(x) := \omega^2 - n^2(1+x^2). $

The following proof is divided into several steps.

{\bf Step 1.}
Taking the real part of the $L^2(-1,1)$-inner product of \eqref{le} with $u_n$ yields
\begin{equation}
\label{li}
 \|u_n'\|^2_{L^2(-1,1)} - \big\langle q_n(x)  u_n,\;u_n\big\rangle_{L^2(-1,1)} = -Re \big\langle F_n,\;u_n\big\rangle_{L^2(-1,1)},
\end{equation}
and taking the imaginary part yields
\begin{equation}
\label{lii}
|\omega|\|u_n\|_{L^2(I_{\eta,\delta})}^2 \le \|F_n\|_{L^2(-1,1)}\|u_n\|_{L^2(-1,1)}.
\end{equation}
If $n > 2|\omega |$, then $-q_n(x) \ge n^2 - \omega^2 >3n^2/4$ for $x\in [-1,1]$, and it follows from \eqref{li} and the Cauchy--Schwarz inequality that
\begin{equation*}
\label{coercive_mid}
    \|u_n'\|_{L^2(-1,1)}^2 + \frac{3}{4}n^2\|u_n\|_{L^2(-1,1)}^2 \le \|F_n\|_{L^2(-1,1)} \|u_n\|_{L^2(-1,1)},
\end{equation*}
which yields  the following estimate
\begin{equation}\label{lle}
    \|u_n'\|_{L^2(-1,1)} + n\|u_n\|_{L^2(-1,1)} \le C\|F_n\|_{L^2(-1,1)},\;\;\;\;n > 2|\omega |.
\end{equation}

{\bf Step 2.}
It remains to consider the case $n \le 2|\omega|$. Set $\kappa:= 1 + |\omega| + n.$ Then
  \begin{equation}
  \label{bq}
  |q_n(x)| \le 2 n^2  + \omega^2   \le   C\kappa^2 \;\; \text{and} \;\;  |q_n'(x)| = |2xn^2| \le C\kappa^2,\;\; \text{for} \; x \in [-1,1].
  \end{equation}
  In the damping region $I_{\eta, \delta}$,   \eqref{le} can be rewritten as
\begin{equation} \label{eq_vn}
     u_n'' + q_n(x)u_n = F_n + i\omega u_n, \quad \text{for } x \in I_{\eta, \delta}.
\end{equation}

Fix a strictly interior subinterval $J \subset I_{\eta, \delta}$, and let  $\phi \in C_0^\infty(I_{\eta, \delta})$ be a  cut-off function such that $0 \le \phi \le 1$ and $\phi \equiv 1$ on $J$.
Multiplying equation \eqref{eq_vn} by $\phi^2 \bar{u}_n$ and integrating by parts over  $I_{\eta, \delta}$, we obtain
\begin{equation}\label{410}
   - \int_{I_{\eta, \delta}} \phi^2 |u_n'|^2 \,dx - 2\int_{I_{\eta, \delta}} \phi \phi' u_n' \bar{u}_n \,dx + \int_{I_{\eta, \delta}} q_n(x) \phi^2 |u_n|^2 \,dx = \int_{I_{\eta, \delta}} (F_n +i\omega u_n)\phi^2 \bar{u}_n \,dx.
\end{equation}
Applying  Young inequality and restricting the  left-hand side integral to $J$, we deduce
\begin{equation}
    \frac{1}{2} \int_{J} |u_n'|^2 \,dx \le \int_{I_{\eta, \delta}} (2|\phi'|^2 + |q_n(x)|\phi^2) |u_n|^2 \,dx + \int_{I_{\eta, \delta}} |F_n + i\omega u_n| |u_n| \phi^2 \,dx.
\end{equation}
Combining this with \eqref{bq} yields
\begin{equation} \label{lej}
    \|u_n'\|_{L^2(J)}^2 \le C \left(\kappa^2  \|u_n\|_{L^2(I_{\eta, \delta})}^2 + \|F_n\|_{L^2(I_{\eta, \delta})}^2 \right).
\end{equation}

For any real-valued $C^1$ function $h$ on an interval $I \subset I_{\eta, \delta}$, multiplying \eqref{eq_vn} by $h \bar{u}_n'$  and integrating by parts over $I$, we obtain
\begin{equation} \label{e_rellich_new}
    2\operatorname{Re} \int_{I} (F_n + i\omega u_n) h \bar{u}_n' \,dx = \left[ h |u_n'|^2 + h q_n |u_n|^2 \right]_{\partial I} - \int_{I} \left( h'|u_n'|^2 + (h q_n)'|u_n|^2 \right) \,dx.
\end{equation}
Fix $x_0 \in J$ and choose $x_1 \in J$ with $x_1 > x_0$. Let $\xi \in C^1([x_0, x_1])$ be a cut-off function satisfying $0 \le \xi \le 1$, with $\xi(x_0) = 1$ and $\xi(x_1) = 0$.
Applying \eqref{e_rellich_new} on $I = [x_0, x_1]$ with $h = \xi$  and using \eqref{bq} yields
\begin{align} \label{lbd}
\begin{split}
    |u_n'(x_0)|^2 + q_n(x_0) |u_n(x_0)|^2 & = -\int_{x_0}^{x_1} \left( \xi' |u_n'|^2 + (\xi q_n)' |u_n|^2 + 2\operatorname{Re}((F_n + i\omega u_n) \xi \bar{u}_n') \right) dx \\
  &  \le   C \int_J \left( |u_n'|^2 + \kappa^2 |u_n|^2 + |F_n|^2 \right) dx.
\end{split}\end{align}
Moreover, since
$$
 |u_n(x_0)|^2 = -\int_{x_0}^{x_1} (\xi |u_n|^2)' \,dx = -\int_{x_0}^{x_1} \left( \xi' |u_n|^2 + 2\xi Re (u_n \bar{u}_n') \right) dx,
$$
we have
 \begin{equation}
\label{lb}
|u_n(x_0)|^2  \le C \int_J \left( |u_n|^2 +  |u'_n| |u_n|  \right) dx.
\end{equation}
By \eqref{lej} and \eqref{lb},
\begin{equation}
\label{lb1}
\kappa^2|u_n(x_0)|^2 \le  C\big(\kappa^3\|u_n\|_{L^2(I_{\eta,\delta})}^2 + \kappa^2\|u_n\|_{L^2(I_{\eta,\delta})}\|F_n\|_{L^2(I_{\eta,\delta})} \big).
\end{equation}
Combining \eqref{lej}, \eqref{lbd} and \eqref{lb1} yields
\begin{equation}\label{lbb}
   \kappa^2|u_n(x_0)|^2 + |u_n'(x_0)|^2 \le C\big(\kappa^3\|u_n\|_{L^2(I_{\eta,\delta})}^2 + \kappa^2\|u_n\|_{L^2(I_{\eta,\delta})}\|F_n\|_{L^2(I_{\eta,\delta})}+\|F_n\|^2_{L^2(I_{\eta,\delta})} \big),\;\; \forall x_0 \in J.
\end{equation}

{\bf Step 3.}
 Let $\mathcal V_{n}(x) :=  |u_{n}'(x)| ^2 + \kappa^2 |u_{n}(x) |^2.$
Then,
\begin{equation} \label{lvp}
\begin{aligned}
    \mathcal{V}_n'(x) &= 2\operatorname{Re}\big( u_n'' \bar{u}_n' \big)  + 2\kappa^2 \operatorname{Re}\big( u_n' \bar{u}_n \big)  \\
    &= 2\operatorname{Re}\big( (- q_n(x) u_n + F_n + i\omega \chi_{I_{\eta, \delta}}u_n) \bar{u}_n' \big)  + 2\kappa^2 \operatorname{Re}\big( u_n' \bar{u}_n \big).
\end{aligned}
\end{equation}
Substituting   \eqref{bq} into \eqref{lvp} and using the bounds $|q_n(x)| \le C\kappa^2$ and $|\omega| \le \kappa$ yields
\begin{equation} \label{lvpp}
\begin{aligned}
   |\mathcal{V}_n'(x)| &\le 2|q_n(x)| |u_n||u_n'| + 2|F_n(x)||u_n'| + 2|\omega| |u_n||u_n'| + 2\kappa^2 |u_n||u_n'| \\
    &\le \tilde C\kappa \mathcal{V}_n(x) + |F_n(x)|^2.
\end{aligned}
\end{equation}
 Therefore,  applying Gronwall's inequality yields
\begin{equation}
    \mathcal{V}_n(x) \le e^{\tilde C\kappa} \left( \mathcal{V}_n(x_0) + \int_{-1}^1 |F_n(s)|^2 \,ds \right), \quad \forall x \in [-1,1].
\end{equation}
Integrating this estimate over $x \in [-1,1]$ and using   \eqref{lbb}, we deduce
\begin{equation} \label{final_gronwall}
    \|u_n'\|_{L^2(-1,1)}^2 + \kappa^2 \|u_n\|_{L^2(-1,1)}^2 \le C e^{\tilde  C\kappa} \left( \kappa^4 \|u_n\|_{L^2(I_{\eta,\delta})}^2 + \|F_n\|_{L^2(-1,1)}^2 \right).
\end{equation}
Since $n \le 2|\omega|$ implies $\kappa \le 1 + 3|\omega|$,  it follows from \eqref{lii}
and \eqref{final_gronwall} that for   $|\omega| \ge 1$,
\begin{equation}
\begin{split}
    \|u_n'\|_{L^2(-1,1)}^2 + \kappa^2 \|u_n\|_{L^2(-1,1)}^2 & \le C e^{\tilde  C\kappa} \left( \kappa^3 \|F_n\|_{L^2(-1,1)}\|u_n\|_{L^2(-1,1)} + \|F_n\|_{L^2(-1,1)}^2 \right)
    \\& \le \frac{1}{2} \kappa^2 \|u_n\|_{L^2(-1,1)}^2 + C  e^{2\tilde  C\kappa} \|F_n\|_{L^2(-1,1)}^2.
\end{split}\end{equation}
Hence, for all sufficiently large $|\omega|$, we obtain the desired estimate
\begin{equation}
    \|u_n'\|_{L^2(-1,1)}^2 + (1+| \omega|+n)^2\|u_n\|_{L^2(-1,1)}^2 \le C e^{\tilde C|\omega|}\|F_n\|_{L^2(-1,1)}^2.
\end{equation}
The case of bounded $|\omega|$
 follows from the continuity of the resolvent on compact sets.
 Finally, summing these mode estimates over $n$, we arrive at \eqref{lr}.
\end{proof}

Lemma \ref{lemma-log} shows that the semigroup $e^{tA} $ is logarithmically stable. The next lemma provides a matching lower bound for the resolvent, which, together with the upper bound in Lemma \ref{lemma-log}, shows that the logarithmic decay rate $\frac{C}{\log(2+t)}$ is optimal. This completes the proof of Theorem \ref{th2.2}.

 \begin{lemma}\label{lemma-log1}
For $n\in {\mathbb N}, \; (x,y) \in   [-1,1] \times [0,\pi]$, let
  \begin{align}\label{dlu}
  \begin{split}
   \omega_n=(n^2+n)^{{1\over2}},  \quad
  \varphi_n(x) =  {\chi(x)\, n^{1 \over 4}e^{-{nx^2 \over 2}} \over \|\chi(x)\, n^{1 \over 4}e^{-{nx^2 \over 2}}\|_{L^2(-1,1)} },  \quad
 U_n(x,y)=\varphi_n(x)e_n(y),
\end{split} \end{align}
where the  cut-off function  $ \chi \in C^2_c[-1,1] $ satisfies $
\|\chi\|_{C[-1,1]} \le 1$ and $\chi \equiv 1 $ on $(-\rho_0,\;\rho_0)$ for some $\rho_0$ with $  \eta-\delta <\rho_0<1.$
   Then
for some constant $C$, independent of $n$,  the normalized vector
$$Z_n={(U_n,\;i\omega_nU_n)\over \|(U_n,\;i\omega_nU_n)\|_{\mathcal H}}$$
satisfies $Z_n \in \mathcal D(A),$ and
\begin{align}\label{lee}
\|(i\omega_n-A)Z_n\|_{\mathcal H}\le  C n^{{1\over4}} e^{-{(\eta-\delta)^2n\over 2}}.
\end{align}
   Consequently,   for suitable constants $C,\;
   \tilde C>0$,   $$\|(i\omega_n-A)^{-1}\|_{\mathcal L(\mathcal H)}\ge C e^{\tilde C(\eta-\delta)^2\omega_n}.$$
 \end{lemma}

\begin{proof}
For the denominator of $Z_n$, we have
\begin{equation}\label{uz}
\begin{split}
      \|(U_n, i\omega_n U_n)\|_{\mathcal{H}}^2 & = \|\partial_x U_n\| ^2 + \|\sqrt{1+x^2}\,\partial_y U_n\| ^2 + \omega_n^2\|U_n\| ^2
      \\ &
      = \|\varphi_n'\|_{L^2(-1,1)}^2 + n^2\|\sqrt{1+x^2}\, \varphi_n\|_{L^2(-1,1)}^2 + \omega_n^2\|\varphi_n\|_{L^2(-1,1)}^2.
      \end{split}
\end{equation}
Let $\phi_n(x) := n^{1 \over 4}e^{-{nx^2 \over 2}}$ and $c_n := \|\chi \phi_n\|_{L^2(-1,1)}^{-1}$.
A direct computation gives
 $$
 \lim_{n \to \infty} \int_{-1}^1  |\chi(x)|^2 n^{1\over2} e^{-nx^2} \, dx
   =\lim_{n \to \infty} \int_{-\sqrt{n}}^{\sqrt{n}} e^{-t^2} \Big|\chi\Big(\frac{t}{\sqrt{n}}\Big)\Big|^2 \, dt
 = |\chi(0)|^2 \int_{-\infty}^{\infty} e^{-t^2} \, dt = |\chi(0)|^2 \sqrt{\pi}.$$
   Thus, there exist positive constants $C_b,\;C_u$ such that
   \begin{equation}
   \label{cb}
   C_b \le c_n \le C_u, \;\; \text{uniformly for }\; n.
   \end{equation}
By definition, $\|\varphi_n\|_{L^2(-1,1)} = 1$. Furthermore, by standard Gaussian integral estimates and \eqref{cb}, we have
\begin{align}
  &\|x\varphi_n\|_{L^2(-1,1)}^2 \le C \int_{-1}^1 x^2 n^{1 \over 2} e^{-nx^2} dx =
{C\over n}\Big(  -\sqrt{n} e^{-n} +{1\over2}  \int_{-\sqrt{n}}^{\sqrt{n}} e^{-t^2} dt \Big)\le C n^{-1},\\ &  \|\varphi_n'\|_{L^2(-1,1)}^2 \le C\big(\|\chi'\phi_n \|_{L^2(-1,1)}^2+ n^2\|\chi x \varphi_n\|_{L^2(-1,1)}^2\big)
    \le C n. \end{align}
Therefore, substituting these into \eqref{uz} yields
\begin{equation}\label{uu}
     \widetilde{C}_1 n^2 \le \|(U_n, i\omega_n U_n)\|_{\mathcal{H}}^2 \le \widetilde{C}_2 n^2,
     \;\;\;\; 0<\widetilde{C}_1<\widetilde{C}_2.
\end{equation}

In view of
\begin{equation}
\label{order}
(i\omega_n - A)\begin{pmatrix} U_n \\ i\omega_n U_n \end{pmatrix} = \begin{pmatrix}0\\ (-\omega_n^2 + \mathbb L + i\omega_n\chi_{\mathbb G})U_n\end{pmatrix},
\end{equation}
we will establish the following two estimates
 \begin{align} \label{p3.2}
  & \|(\omega_n^2 - \mathbb L )U_n\|  \le Cn^{5 \over 4}  e^{- { \rho_0^2   n\over 2}},
   \\ \label{p3.3}
   &
   \|\chi_{\mathbb G} U_n\|  \le Cn^{1 \over 4}e^{-{(\eta-\delta)^2n \over 2}}.
    \end{align}
Since $\rho_0 > \eta-\delta>0$  and  $ \omega_n=(n^2+n)^{{1\over2}}$,  it follows that
    \begin{align} \label{p3.4}
   \|(\omega_n^2 - \mathbb L - i\omega_n\chi_{\mathbb G})U_n\|  \le Cn^{5 \over 4}e^{-{(\eta-\delta)^2n \over 2}}.
     \end{align}
  Combining \eqref{uu}, \eqref{order} and \eqref{p3.4} yields \eqref{lee}. It remains to verify \eqref{p3.2} and \eqref{p3.3}.

Since $\phi_n$ satisfies the  harmonic oscillator equation
\begin{equation} \label{harmonic_osc}
    \phi_n''(x) + (n - n^2 x^2) \phi_n(x)  = 0,
\end{equation}
we have
   \begin{equation} \label{436}
   \begin{split}
    (\omega_n^2 - \mathbb L)U_n & =   \left( \varphi_n'' + (n - n^2x^2)\varphi_n \right) e_n(y)\\
  & = c_n \left( 2\chi'(x)\phi_n'(x) + \chi''(x)\phi_n(x) \right) e_n(y). \end{split}
\end{equation}
Since the supports of $\chi'$ and $\chi''$ are  contained in $\{x \in [-1,1] : |x| \ge \rho_0\}$, we obtain
\begin{equation} \label{437} |\phi_n(x)| \le n^{1 \over 4}e^{-{n\rho_0^2 \over 2}}, \quad |\phi_n'(x)| = |-nx \phi_n(x)| \le n^{5 \over 4}e^{-{n\rho_0^2 \over 2}}, \quad   \text{for}\;\rho_0\le |x| \le 1.
\end{equation}
Estimate \eqref{p3.2} then follows from \eqref{cb}, \eqref{436}-\eqref{437}.

On the other hand, \eqref{p3.3} follows from \eqref{cb} and
\begin{equation}
    \|\chi_{\mathbb G}U_n\|^2 = \int_{\eta-\delta}^{\eta+\delta} |\varphi_n(x)|^2 dx \le C \int_{\eta-\delta}^{\eta+\delta} n^{1 \over 2} e^{-nx^2} dx \le 2\delta C n^{1 \over 2}   e^{-n(\eta-\delta)^2}.
\end{equation}
The proof is complete.
\end{proof}


\begin{remark}  The line \(x = 0\), where the coefficient \(1 + x^2\) attains its minimum, generates a family of stably trapped rays. For the principal symbol \(\rho = \tau^2 - \xi_1^2 - (1+x^2)\xi_2^2\), the bicharacteristic flow, parameterized by \(t\), satisfies \(dx/dt = \xi_1/\tau\), \(d\xi_1/dt = -x\xi_2^2/\tau\), so that \(d^2x/dt^2 = -(\xi_2^2/\tau^2)x\) is a harmonic oscillator. Rays with small transverse momentum \(\xi_1\) oscillate about \(x = 0\) with small amplitude.

If \(|\eta| > \delta\), small-amplitude rays never meet \(\mathbb G\), and Lemma \ref{lemma-log1} constructs Gaussian quasimodes concentrated on those whose damping is exponentially small.
If \(|\eta| < \delta\), the open strip \(\mathbb G\) contains a neighbourhood of \(\{x = 0\}\) and every ray crosses it, so the geometric control condition (GCC) holds.
The borderline case \(|\eta| = \delta\) is the interesting one: the ray \(\{x = 0\}\) never enters the open set \(G\), so GCC fails, yet every nearby ray spends a fixed fraction of its time in \(\mathbb G\). Theorem \ref{th2.1} shows that this suffices for exponential stability.
The proof relies on the resolvent estimate and  the reflection identity \eqref{hre} for the harmonic oscillator.
\end{remark}

 Since the system \eqref{sys} is invariant under the mapping $x \mapsto -x$, our results extend to negative values of $\eta$ as follows
 \begin{theorem}   Assume that $0 < \delta \le 1/4$ and $-1+\delta \le \eta \le 1-\delta$.
\begin{enumerate}
\item[{\rm (i)}] If $-\delta \le \eta \le \delta$, then the semigroup $e^{tA}$ is exponentially stable.

\item[{\rm (ii)}] If $\delta < \eta \le 1-\delta$ or $-1+\delta \le \eta < -\delta$, then the semigroup  $e^{tA}$ is logarithmically stable, and the decay estimate
$ (\log(2+t))^{-1}  $
is optimal.
\end{enumerate}
\end{theorem}

\section{ Galerkin approximation}
\setcounter{equation}{0} \setcounter{theorem}{0}

In Section 3 and 4, we established the stability properties of the continuous system \eqref{sys} under suitable assumptions on the location of the damping region.  In this section, we construct a spectral Galerkin approximation associated with the pseudo-eigenbasis introduced in \eqref{basis}.
We show that the corresponding finite-dimensional semigroups converge to the original semigroup.
  However, the family of Galerkin approximations is not uniformly exponentially stable, even when the continuous system is exponentially stable

  Let $h>0$ be a truncation parameter tending to zero.  Set
\begin{align}\label{pa}
&\mathbb{V}_h := \text{span}\Big\{\phi_{m,n}\;  :\;  \Big( {m\pi \over 2}\Big)^2   + n^2 \le {1\over h^{2}}\Big\}, \\
&  \mathbb{V}(M,n) := \text{span}\{\phi_{m,n}\;  :\;1\le m\le M \} =\text{span} \{\widetilde{e}_1(x)e_n(y),\, \cdots, \, \widetilde{e}_M(x)e_n(y)\},\\
\label{mhn}
& M_h(n) := \max\Big \{m\; :\; \Big( {m\pi \over 2}\Big)^2 + n^2 \le  {1\over h^{2}}\Big\}, \quad M_h(n) = 0 \text{ if it is empty}.
 \end{align}
It follows that
\begin{align}
\mathbb{V}_h = \bigoplus\limits_{n\,: \, M_h(n) \ge 1} \mathbb{V}(M_h(n),n).
 \end{align}

Let $\mathbb{P}_h : L^2(\Omega) \to \mathbb{V}_h$ be the orthogonal projection. Projecting the system \eqref{sys} onto $\mathbb{V}_h$ yields
\begin{equation}
\label{psys}
\ddot{u}_h(t) + \mathbb{P}_h \mathbb{L} \mathbb{P}_h u_h(t) + \mathbb{P}_h \chi_{\mathbb{G}} \mathbb{P}_h \dot{u}_h(t) = 0, \quad u_h(0) = \mathbb{P}_h u_0,\;\dot u_h(0) = \mathbb{P}_h u_1  \in \mathbb{V}_h,
\end{equation}
where  $u_0,\; u_1$ denote  the initial data.
To simplify the notation, we introduce the discrete operators $\mathbb{L}_h := \mathbb{P}_h \mathbb{L} \mathbb{P}_h$ and $B_h := \mathbb{P}_h \chi_{\mathbb{G}} \mathbb{P}_h$    on $\mathbb{V}_h$. Let $U_h(t) = (u_h(t), \dot{u}_h(t))^\top$. Then the corresponding first-order system reads as follows:
\begin{equation}
\label{dabs}
\dot{U}_h(t) = \mathbb{A}_h U_h(t), \quad U_h(0) = (\mathbb{P}_h u_0, \mathbb{P}_h u_1)^\top \in \mathbb{V}_h \times \mathbb{V}_h,
\end{equation}
where $\mathbb{A}_h$ is the finite dimensional generator on $ \mathbb{V}_h \times \mathbb{V}_h$ defined by:
$$ \mathbb{A}_h = \begin{pmatrix} 0 & I \\ -\mathbb{L}_h & -B_h \end{pmatrix}.
$$
The energy space associated with the discrete system \eqref{dabs} is $\mathcal{H}_h = \mathbb{V}_h \times \mathbb{V}_h$, endowed with the norm
$$ \|(u,v)\|_{\mathcal{H}_h}^2 = \langle \mathbb{L}_h u, u \rangle_{\mathbb{V}_h} + \|v\|_{\mathbb{V}_h}^2, \quad \forall (u,v) \in \mathcal{H}_h. $$

For each fixed $n$,~ let
$
u(x,y)=f(x)e_n(y)  \in \mathbb V(M_h(n),n).$
Then,~ $\mathbb Lu=\big(
-f''(x)+n^2(1+x^2)f(x)
 \big)e_n(y) \in \overline{\text{span}}\{\widetilde{e}_m e_n \,:\, m\ge1\},$ which shows that
the operator
\(
 \mathbb L_h
\)
maps
\(
\mathbb V(M_h(n),n)
\)
into itself.
Similarly, since $\chi_{\mathbb G}$ depends only on $x$, the damping operator
$B_h$ also leaves
$\mathbb{V}(M_h(n),n)$ invariant.
Therefore each subspace
$\mathbb V(M_h(n),n)\times \mathbb V(M_h(n),n)$
is invariant under $\mathbb A_h$, which implies the operator $\mathbb A_h$ admits a block diagonal decomposition.

We now turn to the $x$-variable discretization.  For   fixed $h >0 $ and $n\in {\mathbb N}$, let
\begin{align}\label{emkn}
E_{M_h(n) } := \text{span}\{\widetilde{e}_1,\ldots,\widetilde{e}_{M_h(n) }\},
 \end{align}
 and $\mathbb Q_{M_h(n) }$ be the orthogonal projection from $L^2(-1,1)$ onto $E_{M_h(n) }$.
Clearly,  $\mathbb{V}(M_h(n),n)= E_{M_h(n) }\otimes \text{span} \{e_n\}$.

 For the given $h>0, \; n\in {\mathbb N}$, and $M_h(n)$ defined in \eqref{mhn}, we introduce
\begin{align}\label{smhn}
\begin{split}
&D_{M_h(n)} := \text{diag}\Big\{\Big( {\pi \over 2}\Big)^2,\ldots,\Big( {M_h(n)\pi \over 2}\Big)^2\Big\}, \qquad  X_{M_h(n)} := \mathbb  Q_{M_h(n)} x^2 \mathbb  Q_{M_h(n)},
\\&
H_{M_h(n),n} := D_{M_h(n)} + n^2(I + X_{M_h(n)}),  \qquad   B_{M_h(n)} :=\mathbb Q_{M_h(n)} \chi_{(\eta-\delta, \eta+\delta)} \mathbb Q_{M_h(n)} ,
\end{split}
 \end{align}
where $H_{M_h(n),n}  $ and $ B_{M_h(n)}$ represent the conservative and damping matrices, respectively, and
\begin{align}
\big(X_{M_h(n)}\big)_{i,j} = \int_{-1}^1 x^2 \widetilde{e}_i \widetilde{e}_jdx,\;\;\;\;
\big(B_{M_h(n)}\big)_{i,j} = \int_{ \eta-\delta}^{\eta+ \delta}  \widetilde{e}_i \widetilde{e}_jdx.
 \end{align}
 Then substituting the expansion $u_h(x,y,t) = \sum_{j=1}^{M_h(n)} q_j(t) \widetilde{e}_j(x) e_n(y)$ into the projected equation \eqref{psys},
 we obtain the following  finite-dimensional modal system:
 \begin{equation}
 \label{sysdd}
 \ddot{q}(t) + H_{M_h(n),n} q(t)  + B_{M_h(n)} \dot{q}(t)= 0,\;\;\;\;\text{where}\;
 q(t) = (q_1(t), \dots, q_{M_h(n)}(t))^\top .
 \end{equation}

For each pair $h,\;n$, define
$$\widetilde A_{M_h(n),n} = \begin{pmatrix} 0 & I \\ -H_{M_h(n),n} & -B_{M_h(n)} \end{pmatrix},  \qquad
D(\widetilde A_{M_h(n),n}) =\mathbb{V}(M_h(n),n)\times\mathbb{V}(M_h(n),n) .$$
Then the matrix first-order ordinary differential  system corresponding to \eqref{sysdd} is
$$
\dot{Q} (t) = \widetilde A_{M_h(n),n}Q(t), \quad  Q(0) \in  E_{M_h(n)}\times E_{M_h(n)}.
 $$
Let
\begin{align} \label{dlam}
\Lambda_{M_h(n),n}(\omega) = H_{M_h(n),n}   - \omega^2 I   + i \omega B_{M_h(n)}, \;\;\; D(\Lambda_{M_h(n),n}(\omega)) = \mathbb{V}(M_h(n),n).
 \end{align}
 It is clear
that
$i\omega$
is an eigenvalue of
$\widetilde A_{M_h(n),n}$
if and only if
\[
\det\Lambda_{M_h(n),n}(\omega)=0.
\]
Hence, the spectral properties of
$\widetilde A_{M_h(n),n}$
reduce  to the study of the matrix polynomial
$\Lambda_{M_h(n),n}$.

With the above definitions,  $\mathbb A_h$ inherits the block diagonal structure
\begin{equation}
\label{ahd}
\mathbb A_h
=
diag\big\{
\widetilde A_{M_h(1),1}, \,\cdots,\, \widetilde A_{M_h(N),N} \big\},\;\;\text{where } N = \max\{n\,:\,   M_h(n)\ge1\}.
\end{equation}
Consequently, the spectrum of $\mathbb A_h$ is the union of the spectra of the individual blocks:
\[
\sigma(\mathbb A_h)
=
\bigcup_{\,n:\,M_h(n)\ge1}
\sigma\!\left(\widetilde A_{M_h(n),n}\right).
\]



\begin{remark}
Although $\mathbb{L}_0$ commutes with the projection,
$\mathbb{P}_h\mathbb{L}_0 = \mathbb{L}_0\mathbb{P}_h$, the variable coefficient
operator $\mathbb{L}$ does \emph{not} commute with $\mathbb{P}_h$. For instance,
since $\mathbb{L} = \mathbb{L}_0 - x^2\partial_{yy}$ and
$\langle x^2\widetilde{e}_1,\widetilde{e}_3\rangle_{L^2(-1,1)} = 3/(2\pi^2)\neq 0$,
one cannot expect $\mathbb{P}_h\mathbb{L} = \mathbb{L}\mathbb{P}_h$.
\end{remark}



\begin{theorem} For every $(u, \, v)\in \mathcal{H}$ and every $T>0$, it holds
\begin{equation}
\label{dc}
\lim_{h\to0}\,
\sup_{0\le t\le T} \big\| e^{t\mathbb{A}_h}  \mathcal{P}_h(u,\,v) - e^{tA} (u, \, v) \big\|_{\mathcal H}
=0.  \end{equation}
\end{theorem}

\begin{proof}
{\bf (i)} Set $
\mathcal{P}_h(u,\, v) := (\mathbb{P}_h u,\, \mathbb{P}_hv) $ for$  (u,\,v) \in \mathcal{H}.$
Since $
\{\phi_{m,n}\}
$
forms a complete orthonormal basis of $L^2(\Omega)$, we have by the same argument as Remark \ref{remark-ll},
\[
\big\| \mathbb{L}^{1\over2}  \mathbb{P}_h u\|^2
 \le 2 \big\| \mathbb{L}^{1\over2}_0 \mathbb{P}_h u\|^2
 \le 2 \big\| \mathbb{L}^{1\over2}_0   u\|^2  \le 2 \big\| \mathbb{L}^{1\over2}    u\|^2  , \;\;\; \forall \; u\in H^1_0(\Omega).
\]
Consequently, $\|\mathcal{P}_h(u,\, v)\|_{\mathcal H} \le C \| (u,\, v)\|_{\mathcal H} $ and
\begin{align}\label{reso1}
 \lim_{h \to 0} \|\mathcal{P}_h(u,\, v)  - (u,\,v)\|_{\mathcal{H}} = 0,\;\; \forall \;(u,\,v) \in \mathcal{H}.
 \end{align}

{\bf  (ii)} Let $\varepsilon_h(u,v) := \mathbb{A}_h \mathcal{P}_h(u,v) - \mathcal{P}_h A(u,v)$   for  $(u,\,v) \in D(A)$.
  It follows that
\[
\varepsilon_h(u,v) = \big(0,\;
\mathbb{P}_h \mathbb{L}(u - \mathbb{P}_h u) + \mathbb{P}_h \chi_{\mathbb G} (v - \mathbb{P}_h v) \big).
\]
Since $u \in D(\mathbb{L})$,  one has  $\sum_{m,n}\lambda_{m,n}^2|u_{m,n}|^2<\infty$, where $u_{m,n}$ is defined in $u = \sum u_{m,n} \phi_{m,n}$.
 Therefore,
\begin{equation}
\label{re1}\|\mathbb{L}(u - \mathbb{P}_h u)\|^2 \le 4 \|\mathbb{L}_0(u- \mathbb{P}_h u)\|^2
= 4 \sum_{ \lambda_{m,n}>h^{-2}} \lambda_{m,n}^2 |u_{m,n}|^2   \to 0.
\end{equation}
Similarly,
\begin{equation}
\label{re2}
\|\chi_{\mathbb G} (v - \mathbb{P}_h v)\| \le    \|v - \mathbb{P}_h v\|  \to 0.
\end{equation}
 Hence, by  \eqref{re1}, \eqref{re2}, we obtain as $h\to0,$
  \begin{equation}
  \label{re}
   \varepsilon_h(u,v) \to 0  \;  \mbox{   in}  \; \mathcal{H}.
  \end{equation}

 {\bf  (iii)} In this step, we shall prove that for all $U \in \mathcal{H} $ and $\omega >0 ,$
 \begin{align} \label{reso}
\lim_{h \to 0} \|\big(R(\omega  I,\, \mathbb A_h)\mathcal{P}_h- R(\omega  I,\,  A) \big) U\|_{\mathcal H }
= 0.
\end{align}
Since
\[
(\omega I - \mathbb{A}_h) \mathcal{P}_h  R(\omega, A) U
= \mathcal{P}_h(\omega  R(\omega, A) U - A R(\omega, A) U) - \varepsilon_h( R(\omega, A) U)
= \mathcal{P}_hU  - \varepsilon_h( R(\omega, A) U),
\]
we have
\begin{equation}
\label{re3}
R(\omega, \mathbb{A}_h) \mathcal{P}_hU -  R(\omega, A) U
= \big(\mathcal{P}_h  R(\omega, A) U -  R(\omega, A) U\big)+ R(\omega, \mathbb{A}_h) \varepsilon_h( R(\omega, A) U).
\end{equation}
Then, by \eqref{reso1}, \eqref{re}, \eqref{re3} and the fact that $\mathbb A_h$ generates a
$C_0$-semigroup of contractions, we obtain
\begin{equation}
\label{re4}
\|R(\omega, \mathbb{A}_h) \mathcal{P}_hU - R(\omega, A)U\|_{\mathcal{H}}
\le \|\mathcal{P}_h  R(\omega, A) U -  R(\omega, A) U\|_{\mathcal{H}} + \frac{1}{\omega} \|\varepsilon_h( R(\omega, A) U)\|_{\mathcal{H}}
\to 0 .
\end{equation}
 Hence, \eqref{reso} is proved.

{\bf (iv)}
 For $U\in\mathcal{H}$ and $T>0$,
we decompose
\begin{equation}\label{d123}
\big(e^{t\mathbb{A}_h}\mathcal{P}_h-e^{tA}\big)R(\omega,A)U=D_1(t)+D_2(t)+D_3(t),
\end{equation}
where
\[
\begin{aligned}
D_1(t)&=e^{t\mathbb{A}_h}\mathcal{P}_h\big(R(\omega,A)-R(\omega,\mathbb{A}_h)\mathcal{P}_h\big)U,\\
D_2(t)&=R(\omega,\mathbb{A}_h)\big(e^{t\mathbb{A}_h}\mathcal{P}_h-\mathcal{P}_he^{tA}\big)U,\\
D_3(t)&=\big(R(\omega,\mathbb{A}_h)\mathcal{P}_h-R(\omega,A)\big)e^{tA}U.
\end{aligned}
\]
Due to Remark \ref{remark-ll}, \eqref{reso} and the fact that $\{e^{tA}U:0\le t\le T\}$ is compact in $\mathcal{H}$, we have
\begin{align}\label{d124}
\begin{split}
&\sup_{0\le t\le T}\|D_1(t)\|_{\mathcal{H}}\le
\sqrt{2}\,\big\|\big(R(\omega,A)-R(\omega,\mathbb{A}_h)\mathcal{P}_h\big)U\big\|_{\mathcal{H}}\to0,\\
&
\sup_{0\le t\le T}\|D_3(t)\|_{\mathcal{H}}\to0.
\end{split}
 \end{align}

To estimate $D_2$, we   introduce the function
\[
f_h(s):=-\,e^{(t-s)\mathbb{A}_h}R(\omega,\mathbb{A}_h)\,\mathcal{P}_h\,e^{sA}R(\omega,A)U,\;\;\;\;
U\in\mathcal{H},\; 0\le s\le t,
\]
which is differentiable.
A direct computation gives
\[
\frac{d}{ds}f_h(s)
=e^{(t-s)\mathbb{A}_h}\big[R(\omega,\mathbb{A}_h)\mathcal{P}_h-\mathcal{P}_hR(\omega,A)\big]e^{sA}U.
\]
Integrating over $s\in[0,t]$ and using that $\mathcal{P}_h$ acts as the identity
on the range of $R(\omega,\mathbb{A}_h)$, we obtain
\begin{equation}\label{duhamel}
D_2(t)R(\omega,A)U
=\int_0^t e^{(t-s)\mathbb{A}_h}\big[R(\omega,\mathbb{A}_h)\mathcal{P}_h-\mathcal{P}_hR(\omega,A)\big]e^{sA}U\,ds.
\end{equation}
By \eqref{reso1} and \eqref{reso},
\[
R(\omega,\mathbb{A}_h)\mathcal{P}_h-\mathcal{P}_hR(\omega,A)
=\big[R(\omega,\mathbb{A}_h)\mathcal{P}_h-R(\omega,A)\big]+(I-\mathcal{P}_h)R(\omega,A)\to0,
\quad\text{strongly}.
\]
Combining this with the boundedness of the integrand, we can deduce from  the Lebesgue dominated convergence theorem   that, for each $t\in[0,T]$,
\begin{equation}\label{d2conv}
\lim_{h\to0}\big\|D_2(t)R(\omega,A)U
\big\|_{\mathcal{H}}=0,
\qquad \forall\,U\in\mathcal{H}.
\end{equation}
Consequently, by \eqref{d124} and \eqref{d2conv},
\begin{equation}\label{conv-d2}
\lim_{h\to0}\sup_{0\le t\le T}\big\|\big(e^{t\mathbb{A}_h}\mathcal{P}_h-e^{tA}\big)V\big\|_{\mathcal{H}}=0,
\qquad \forall\,V\in D(A^2).
\end{equation}
Finally, since $D(A^2)$ is dense in
$\mathcal{H}$,  a standard  argument yields \eqref{dc} for every $(u, v) \in H$.
\end{proof}

To study the uniform exponential stability of the approximation semigroups $\{e^{t\mathbb{A}_h}\}_{h>0}$, we use a  recent result from \cite{lks-l}, which establishes the equivalence between the uniform exponential stability of modal approximation systems and a uniform exact observability inequality.

\begin{lemma}
\label{lemma-g3}
The family of systems \eqref{psys} is uniformly exponentially stable
if and only if  there exists a constant \(C > 0\) such that
\begin{equation}\label{526}
\|(\mathbb L_h - \omega^2)u_h\|  + \|\omega \chi_{\mathbb G} u_h\|  \geq C \|\omega u_h\|
\end{equation}
 for all $\omega \in \mathbb{R}, u_h \in \mathbb V_h $ and $ h>0$.
\end{lemma}

The following result   shows that the family of semigroups $\{e^{t\mathbb{A}_h}\}_{h>0}$  is not uniformly exponentially stable.

 \begin{theorem}\label{th-ne} Assume $|\eta|+\delta<1. $ The family  of   semigroups $\{e^{t\mathbb{A}_h}\}_{h>0}$ is not uniformly exponentially stable.
 \end{theorem}

 \begin{proof}  For each $M \in {\mathbb N}$, since $X_M$ defined in \eqref{smhn} is a self-adjoint operator on  finite-dimensional  space $E_M$, denote by
  $\xi_M := \max\sigma(X_M)  $  its largest eigenvalue, and let
 \begin{equation}
 \label{xq}X_M q_M = \xi_M q_M, \;\; \|q_M\|_{L^2(-1,1)} = 1.
 \end{equation}
 Then,
 \begin{equation}\label{qm}
\begin{aligned}
\xi_M
&=
\max_{\substack{u_M\in E_M\\ \|u_M\|_{L^2(-1,1)}=1}}
\big\langle X_Mu_M,u_M\big\rangle_{L^2(-1,1)}    =
\max_{\substack{u_M\in E_M\\ \|u_M\|_{L^2(-1,1)}=1}}
\int_{-1}^1x^2|u_M(x)|^2\,dx .
\end{aligned}
\end{equation}
It follows that
\[
0\le\xi_M\le1.
\]
Moreover, the inclusion $E_M\subset E_{M+1}$ implies that the sequence
$\{\xi_M\}_{M\ge1}$ is nondecreasing.
 Introduce
\[
\varphi_\varepsilon\in C_c^\infty(-1,1),
\qquad
\operatorname{supp}\varphi_\varepsilon
\subset(1-\varepsilon,1),
\qquad
\|\varphi_\varepsilon\|_{L^2(-1,1)}=1,  \qquad  \varepsilon>0.
\]
and
\[
u_M
=
\frac{\mathbb Q_M\varphi_\varepsilon}
{\|\mathbb Q_M\varphi_\varepsilon\|_{L^2(-1,1)}}
\in E_M.
\]
It follows that as $M\to \infty,$
\[
u_M
\to
\varphi_\varepsilon
\qquad
\text{strongly in }L^2(-1,1).
\]
Let  $ T f(x) = x^2 f(x)$. Since $T$ is bounded,
\[
\big\langle Tu_M, u_M \big\rangle_{L^2(-1,1)} \to \big\langle T \varphi_\varepsilon, \varphi_\varepsilon \big\rangle_{L^2(-1,1)}.
\]
Furthermore,
\[
\big\langle T\varphi_\varepsilon,\varphi_\varepsilon\big\rangle_{L^2(-1,1)}
=
\int_{1-\varepsilon}^1x^2|\varphi_\varepsilon|^2dx
\ge
(1-\varepsilon)^2.
\]
Combining this with \eqref{qm}, we obtain
$
\liminf_{M \to \infty} \xi_M \geq (1-\varepsilon)^2.
$
Since $\varepsilon>0$ is arbitrary and $\xi_M\le1$, we conclude that
\begin{align}\label{xi}
\xi_M \to 1.
 \end{align}

  Next define
$$d_M := \big\langle B_M q_M, \, q_M \big\rangle_{L^2(-1,1)} =\big\langle \mathbb Q_{M} \chi_{(\eta-\delta, \eta+\delta)} \mathbb Q_{M}q_M, \, q_M \big\rangle_{L^2(-1,1)}.
 $$
 Then, similar to \eqref{qm},
 we have
 \begin{align}\label{dmd}
 d_M = \int_{\eta-\delta}^{\eta+\delta} |q_M(x)|^2\,dx.
 \end{align}
Moreover, by \eqref{qm} and \eqref{dmd},
$$(1-(|\eta|+\delta)^2)d_M \le \int_{-1}^1 (1-x^2)|q_M(x)|^2\,dx \le 1-\xi_M .$$
Combining this  with \eqref{xi}
yields
\begin{align}
\label{dm0}d_M \le \frac{1-\xi_M}{ 1-(|\eta|+\delta)^2 } \to 0.
 \end{align}

Now   for each $M$, let
\begin{equation}
\label{nm}n_M=M^3,
\end{equation}
and choose $h_M>0$ such that
\begin{equation}
\label{nm1}
\left(\frac{M\pi}{2}\right)^2+n_M^2
\le
h_M^{-2}
<
\left(\frac{(M+1)\pi}{2}\right)^2+n_M^2.
\end{equation}
Then, by the definition of $M_h(n)$ in \eqref{mhn},
\begin{equation}
\label{nm2}
M_{h_M}(n_M)=M.
\end{equation}

Define \begin{align}\label{hq}\omega_M^2 := \big\langle H_{M,n_M} q_M, q_M \big\rangle_{L^2(-1,1)}= \big\langle D_M q_M, q_M \big\rangle_{L^2(-1,1)} + n_M^2(1+\xi_M).
 \end{align}
It follows that
\begin{align} \label{hm}\begin{split}
  (H_{M,n_M} - \omega_M^2 I)q_M
=  D_M  q_M - \big\langle D_M q_M, q_M \big\rangle_{L^2(-1,1)}   q_M.
\end{split} \end{align}
Since $\|D_M\|= \big(M\pi/2\big)^2 \le CM^2$ and $\omega_M \ge  \big\langle H_{M,n_M} q_M, q_M \big\rangle_{L^2(-1,1)}^{1\over2} \ge n_M = M^3$, we can conclude from \eqref{xq} and \eqref{hm} that
\begin{align}\label{rm}
{\|(H_{M,n_M} - \omega_M^2 I)q_M\|_{L^2(-1,1)} \over \omega_M}\le
{CM^2 \over M^3} \to 0.
\end{align}
Furthermore, due to \eqref{dmd}, \eqref{dm0}, we have
\begin{equation}
\label{bm0}\|B_Mq_M\|_{L^2(-1,1)}^2
\le
\big\langle B_Mq_M,q_M\big\rangle_{L^2(-1,1)}
=
d_M \to 0 .
\end{equation}

Finally, we embed $q_M/\omega_M$ into the full Galerkin space. To this end, let
\[
y_M(x,y) := \frac{1}{\omega_M} q_M(x) e_{n_M}(y) \in \mathbb{V}(M, n_M) \subset \mathbb{V}_{h_M}.
\]
 Then, using \eqref{rm} and \eqref{bm0} yields
\begin{align}
&\| (\mathbb L_{h_M} - \omega_M^2) y_M \|  + \| \omega_M B_{h_M} y_M \|  \to0,\;\; \text{as}\; M\to\infty,\\
& \| \omega_M y_M \| =\| q_M \|_{L^2(-1,1)} \| e_{n_M} \|_{L^2(0,\pi)} =1.
 \end{align}
This contradicts  \eqref{526} for large $M$, which proves that the family of Galerkin semigroups $\{e^{t\mathbb{A}_h}\}_{h>0}$ is not uniformly exponentially stable.
\end{proof}

\begin{remark}
The loss of uniform exponential stability in Theorem~\ref{th-ne} is caused
by the choice of the basis \eqref{basis}. Let $\{\Phi_k\}_{k\ge1}$ be an
orthonormal basis of $L^2(\Omega)$ consisting of eigenfunctions of
$\mathbb L$ itself, with $\mathbb L\Phi_k=\mu_k\Phi_k$. Separation of
variables gives $\Phi_k=\psi_{j,n}(x)e_n(y)$, where the $\psi_{j,n}$ are the
Dirichlet eigenfunctions of $-\frac{d^2}{dx^2}+n^2(1+x^2)$ on $(-1,1)$. Let
$\Pi_N$ denote the orthogonal projection onto
$\operatorname{span}\{\Phi_1,\dots,\Phi_N\}$. Then $\Pi_N$ commutes with
$\mathbb L$, so that
$\Pi_N\mathbb L\Pi_N=\mathbb L|_{\operatorname{Ran}\Pi_N}$: the variable
coefficient, i.e., multiplication by $x^2$ in
$\mathbb L=\mathbb L_0-x^2\partial_{yy}$, is absorbed into the modes and does
not appear in the approximating systems. Only the bounded damping $\chi_G$ is
approximated, through its compression $\Pi_N\chi_G\Pi_N$. In this setting, it
follows from \cite[Theorem~3.6]{lks-l} that the modal approximation is
uniformly exponentially stable if and only if the continuous system is
exponentially stable; under the assumptions of Theorem~4.5, the latter holds
if and only if $|\eta|\le\delta$. In short, the true modes diagonalize the
full operator $\mathbb L$, variable coefficient included, whereas the
pseudo-modes diagonalize only $\mathbb L_0$. The finite-section approximation
of the multiplication by $x^2$ then creates spurious, essentially undamped
modes near the boundary $x=\pm1$.

\end{remark}

\begin{remark}
In Theorem \ref{th-ne}, we have   proved that the family of semigroups fails to be uniformly exponentially stable.
Consequently, an interesting problem for future work is to determine whether the system is uniformly polynomially stable or only logarithmically stable.
Resolving this question will likely require a delicate estimate of the resolvent of $\widetilde A_{M_h(n),n}$ along the imaginary axis, which we leave as a subject for future investigation.
\end{remark}


\bibliographystyle{plain}

\begin{thebibliography}{20}
\bibitem{nicaise}
F.~Abdallah, S.~Nicaise, J.~Valein and A.~Wehbe,
\newblock Uniformly exponentially or polynomially stable approximations for second order evolution equations and some applications,
\newblock {\em ESAIM: COCV},
19(3): 844--887, 2013.

\bibitem{AL14}
N.~Anantharaman and M.~L\'eautaud,
\newblock Sharp polynomial decay rates for the damped wave equation on the torus,
\newblock {\em Anal. PDE},
7(1): 159--214, 2014.



\bibitem{bardos}
C.~Bardos, G.~Lebeau, and J.~Rauch,
\newblock Sharp sufficient conditions for the observation, control, and stabilization of waves from the boundary,
\newblock {\em SIAM J. Control Optim.},
30: 1024--1065, 1992.

\bibitem{batty_duyckaerts}
C.~J.~K.~Batty and T.~Duyckaerts,
\newblock Non-uniform stability for bounded semi-groups on Banach spaces,
\newblock {\em J. Evol. Equ.},
8: 765--780, 2008.

\bibitem{bt}
A.~Borichev and Y.~Tomilov,
\newblock Optimal polynomial decay of functions and operator semigroups,
\newblock {\em Math. Ann.},
347: 455--478, 2010.


\bibitem{Bur98}
N. Burq.
\newblock D\'ecroissance de l'\'energie locale de l'\'equation des ondes pour le probl\`eme ext\'erieur et absence de r\'esonance au voisinage du r\'eel.
\newblock \textit{Acta Math.}, 180(1):1--29, 1998.

\bibitem{BG97}
N. Burq and P. G\'erard.
\newblock Condition n\'ecessaire et suffisante pour la contr\^olabilit\'e exacte des ondes.
\newblock \textit{C. R. Acad. Sci. Paris S\'er. I Math}, 325(7):749--752, 1997.


\bibitem{BH07}
N. Burq and M. Hitrik.
\newblock Energy decay for damped wave equations on partially rectangular domains.
\newblock \textit{Mathematical Research Letters}, 14(1):35--47, 2007.



%
%

\bibitem{chenliuliu}
S.~Chen, K.~Liu, and Z.~Liu,
\newblock Spectrum and stability for elastic systems with global or local Kelvin-Voigt damping,
\newblock {\em SIAM J. Appl. Math.},
59(2): 651--668, 1998.

\bibitem{CSVW14}
H. Christianson, E. Schenck, A. Vasy, and J. Wunsch.
\newblock From resolvent estimates to damped waves.
\newblock \textit{J. Anal. Math.}, 122(1):143--162, 2014.

\bibitem{cox}
S.~Cox and E.~Zuazua,
\newblock The rate at which energy decays in a damped string,
\newblock {\em Commun. Partial Differ. Equ.},
19: 213--243, 1994.

\bibitem{L10}
K.~Datchev and P.~Kleinhenz,
\newblock Sharp polynomial decay rates for the damped wave equation with H\"older-like damping,
\newblock {\em Proc. Amer. Math. Soc.},
148(8): 3417--3425, 2020.

\bibitem{engel}
K.-J.~Engel and R.~Nagel,
\newblock One-parameter semigroups for linear evolution equations,
\newblock {\em Springer},
2000.

\bibitem{ez}
S.~Ervedoza and E.~Zuazua,
\newblock Uniformly exponentially stable approximations for a class of damped systems,
\newblock {\em J. Math. Pures Appl.},
91(1): 20--48, 2009.

\bibitem{fxy}
X.~Fu and L.~Kong,
\newblock Longtime Behavior of Wave Equation with Kinetic Boundary Condition,
\newblock {\em Applied Mathematics \& Optimization},
84: 2803--2817, 2021.

\bibitem{gearhart}
L.~Gearhart,
\newblock Spectral theory for contraction semigroups on Hilbert space,
\newblock {\em Trans. Amer. Math. Soc.},
236: 385--394, 1978.

\bibitem{huang}
F.~Huang,
\newblock Characteristic conditions for exponential stability of linear dynamical systems in Hilbert spaces,
\newblock {\em Ann. Differential Equations},
1: 43--56, 1985.

\bibitem{zua99}
J.~A. Infante and E.~Zuazua,
\newblock Boundary observability for the space semi-discretizations of the 1-D wave equation,
\newblock {\em M2AN},
33(2): 407--438, 1999.

\bibitem{L11}
O.~Kavian and Q.~Zhang,
\newblock Polynomial stability of an abstract system with local damping,
\newblock {\em J. Math. Anal. Appl.},
512(2): 126133, 2022.

\bibitem{L7}
M.~L\'eautaud and N.~Lerner,
\newblock Energy decay for a locally undamped wave equation,
\newblock {\em Ann. Fac. Sci. Toulouse Math.},
(6) 26(1): 157--205, 2017.

\bibitem{Leb96}
G.~Lebeau,
\newblock \'{E}quation des ondes amorties,
\newblock {\em Algebraic and geometric methods in mathematical physics (Kaciveli, 1993), Math. Phys. Stud., 19, Kluwer Acad. Publ., Dordrecht},
73--109, 1996.



\bibitem{L3}
K.~Liu,
\newblock Locally distributed control and damping for the conservative systems,
\newblock {\em SIAM J. Control Optim.},
35: 1574--1590, 1997.

\bibitem{lks-l}
K.~Liu  and L.~H.~ Lu,
\newblock Uniform observability  and uniformly exponential stabilization for a sequence of elastic systems under the model approximation,
submitted.

\bibitem{lzy-zhsm}
Z.~Y. Liu and S.~M. Zheng,
\newblock Uniform exponential stability and approximation in control of a thermoelastic system,
\newblock {\em SIAM J. Control Optim.},
32: 1226--1246, 1994.



\bibitem{lyubich_vu}
Yu.~I.~Lyubich and V.~Q.~Ph\'{o}ng,
\newblock Asymptotic stability of linear differential equations in Banach spaces,
\newblock {\em Studia Math.},
88: 37--42, 1988.

\bibitem{L14}
A.~M\"unch and A.~F.~Pazoto,
\newblock Uniform stabilization of a viscous numerical approximation for a locally damped wave equation,
\newblock {\em ESAIM: COCV},
13(2): 265--293, 2007.


\bibitem{pruss}
J.~Pr\"{u}ss,
\newblock On the spectrum of $C_0$-semigroups,
\newblock {\em Trans. Amer. Math. Soc.},
284: 847--857, 1984.

\bibitem{tucsnak}
K.~Ramdani, T.~Takahashi and M.~Tucsnak,
\newblock Uniformly exponentially stable approximations for a class of second order evolution equations---application to LQR problems,
\newblock {\em ESAIM: COCV},
13(3): 503--527, 2007.


\bibitem{L12}
M.~Rouveyrol,
\newblock Stabilization of the wave equation on larger-dimension tori with rough dampings,
\newblock {\em Pure Appl. Anal.},
6(2): 487--520, 2024.




\bibitem{Ral69}
J. Ralston.
\newblock Solutions of the wave equation with localized energy.
\newblock \textit{Communications on Pure and Applied Mathematics}, 22(6):807--823, 1969.

\bibitem{RT75}
J. Rauch and M. Taylor.
\newblock Exponential decay of solutions to hyperbolic equations in bounded domains.
\newblock \textit{Indiana Univ. Math. J.}, 24(1):79--86, 1974.

\bibitem{L6}
R.~Stahn,
\newblock Optimal decay rate for the wave equation on a square with constant damping on a strip,
\newblock {\em Z. Angew. Math. Phys.},
68(2):  Article 36, 10, 2017.

\bibitem{L13}
L.~R.~Tcheugou\'e T\'ebou and E.~Zuazua,
\newblock Uniform exponential long time decay for the space semi-discretization of a locally damped wave equation via an artificial numerical viscosity,
\newblock {\em Numer. Math.},
95(3): 563--598, 2003.

\bibitem{thang}
S.~Thangavelu,
\newblock Lectures on Hermite and Laguerre expansions,
\newblock {\em Princeton University Press},
1993.


\bibitem{L8}
H.~Zhu,
\newblock Stabilization of damped waves on spheres and Zoll surfaces of revolution,
\newblock {\em ESAIM: COCV},
24(4): 1759--1788, 2018.
%

\bibitem{zua05}
E.~Zuazua,
\newblock Propagation, observation, and control of waves approximated by finite difference methods,
\newblock {\em SIAM Rev.},
47(2): 197--243, 2005.

 \end{thebibliography}

\end{document}